\documentclass{article}
\usepackage{amsthm}
\usepackage{fleqn}
\usepackage{array}
\usepackage[latin1]{inputenc}
\usepackage{amssymb}
\usepackage{amsfonts}
\usepackage{amsmath}

\begin{document}

\newtheorem{thm}{Theorem}[section]
\newtheorem{prop}{Proposition}[section]
\newtheorem{lem}{Lemma}[section]
\newtheorem{cor}{Corollary}[section]

\numberwithin{equation}{section}

\begin{center}
{\Large Rotational surfaces in a normed $3$-space 

\vspace{2mm}

whose principal curvatures satisfy a linear relation }

\vspace{5mm}

Makoto SAKAKI and Kakeru YANASE 
\end{center}

\vspace{2mm}

{\bf Abstract.} We classify rotational surfaces in a normed $3$-space with rotationally symmetric norm whose principal curvatures satisfy a linear relation.

\vspace{2mm}


\vspace{2mm}


\section{Introduction}

It is interesting to generalize differential geometry of surfaces in the Euclidean $3$-space to that in normed $3$-spaces (cf. [1], [2], [3], [4], [8]). In the Euclidean $3$-space, a surface is said to be of linear Weingarten type if its principal curvatures satisfy a linear relation, or if its Gaussian and mean curvatures satisfy a linear relation (cf. [5], [6], [7]). 

In this paper, we define surfaces of linear Weingarten type in a normed $3$-space as surfaces whose principal curvatures satisfy a linear relation, and classify rotational surfaces of linear Weingarten type in a normed $3$-space with rotationally symmetric norm. It is a generalization of [7] in the Euclidean case, and also a generalization of [8] for rotational surfaces with constant mean curvature in a normed $3$-space. 

This paper is organized as follows. In Section 2, basic facts on surfaces in normed $3$-spaces are given (cf. [2]). In Section 3, we give some computations on rotational surfaces in a normed $3$-space with rotationally symmetric norm (cf. [8]). In Section 4, we discuss rotational surfaces in the normed $3$-space one of whose principal curvatures is constant. In Section 5, we consider rotational surfaces of homogenous linear Weingarten type in the normed $3$-space. In Section 6, we study rotational surfaces of inhomogeneous linear Weingarten type in the normed $3$-space.

\section{Preliminaries}

Let $({\mathbb R}^3, \|\cdot\|)$ be a $3$-dimensional normed space. The unit ball $B$ and the unit sphere $S$ are given by
\[B = \{x \in {\mathbb R}^3; \| x\| \leq 1\}, \ \ \ \ S = \{x \in {\mathbb R}^3; \| x\| = 1\}. \]
We assume that $S$ is smooth and strictly convex, that is, $S$ is a smooth surface and $S$ contains no line segment. 

Let $v$ be a non-zero vector in ${\mathbb R}^3$ and $\Pi$ be a plane in ${\mathbb R}^3$. The vector $v$ is said to be Birkhoff orthogonal to $\Pi$, denoted by $v \dashv_{B} \Pi$, if the tangent plane of $S$ at $v/\|v\|$ is parallel to $\Pi$. 

Let $M$ be a surface immersed in $({\mathbb R}^3, \|\cdot\|)$, and $T_{p}M$ be the tangent plane of $M$ at $p \in M$. Then there exists a vector $\eta(p) \in S$ such that $\eta(p) \dashv_{B} T_{p}M$. This gives a local smooth map $\eta: U \subset M \rightarrow S$, which is called the Birkhoff-Gauss map. It can be global if and only if $M$ is orientable. 

Let $k_1$ and $k_2$ be the eigenvalues of $d\eta_p$, which are called the principal curvatures. We say that the surface $M$ is of linear Weingarten type if $a k_1+b k_2 = c$ for some constants $a$, $b$ and $c$ which are not all zero.

\section{Rotational surfaces}

Let $({\mathbb R}^3, \|\cdot\|)$ be a normed $3$-space with rotationally symmetric norm
\[\|x\| = \left( (x_{1}^{2}+x_{2}^{2})^{m}+x_{3}^{2m} \right)^{\frac{1}{2m}}, \ \ \ x = (x_1, x_2, x_3), \]
where $m$ is a positive integer. Set
\[\Phi(x_1, x_2, x_3) := (x_{1}^{2}+x_{2}^{2})^m+x_{3}^{2m}. \]
Then the unit sphere $S$ is given by
\[S = \{(x_1, x_2, x_3) \in {\mathbb R}^3 | \Phi(x_1, x_2, x_3) = 1 \}. \]
Since the case $m = 1$ is the Euclidean case, we assume that $m \geq 2$ in the following. 

Let $M$ be a surface in the above $({\mathbb R}^3, \|\cdot\|)$ which is rotational around $x_3$-axis, and is parametrized by
\[X(u, v) = (\alpha(u)\cos{v}, \alpha(u)\sin{v}, \beta(u)) \]
where $\alpha > 0$, $\alpha' \neq 0$ and $\beta' \neq 0$. Then
\[X_u = (\alpha'\cos{v}, \alpha'\sin{v}, \beta'), \ \ \ \ X_v = (-\alpha\sin{v}, \alpha\cos{v}, 0). \]
The Birkhoff-Gauss map $\eta = \eta(u, v)$ is characterized by the condition
\[(\mbox{grad}(\Phi))_{\eta} = \left( \frac{\partial \Phi}{\partial x_1}(\eta), \frac{\partial \Phi}{\partial x_2}(\eta), \frac{\partial \Phi}{\partial x_3}(\eta) \right) = \varphi X_u \times X_v, \]
where $\varphi$ is a positive function and $\times$ is the standard cross product in ${\mathbb R}^3$. Then we have
\[\eta = A^{-\frac{1}{2m}} \left( -(\beta')^{\frac{1}{2m-1}}\cos{v}, -(\beta')^{\frac{1}{2m-1}}\sin{v}, (\alpha')^{\frac{1}{2m-1}} \right) \]
where
\[A:= \left( \alpha' \right)^{\frac{2m}{2m-1}}+\left( \beta' \right)^{\frac{2m}{2m-1}}. \]

\vspace{1mm}

We can compute that
\[\eta_u = -\frac{1}{2m-1} A^{-\frac{2m+1}{2m}} (\alpha')^{-\frac{2m-2}{2m-1}} (\beta')^{-\frac{2m-2}{2m-1}} (\alpha'\beta''-\alpha''\beta') X_u \]
and
\[\eta_v = -\frac{1}{\alpha} A^{-\frac{1}{2m}}(\beta')^{\frac{1}{2m-1}} X_v, \]
so that
\[k_1 = -\frac{1}{2m-1} A^{-\frac{2m+1}{2m}} (\alpha')^{-\frac{2m-2}{2m-1}} (\beta')^{-\frac{2m-2}{2m-1}} (\alpha'\beta''-\alpha''\beta') \]
and
\[k_2 = -\frac{1}{\alpha} A^{-\frac{1}{2m}}(\beta')^{\frac{1}{2m-1}}. \]
In particular, when $\beta(u) = u$, we have
\begin{eqnarray}
k_1 = \frac{1}{2m-1} \left( \left( \alpha' \right)^{\frac{2m}{2m-1}}+1 \right)^{-\frac{2m+1}{2m}} (\alpha')^{-\frac{2m-2}{2m-1}} \alpha'', 
\end{eqnarray}
and
\begin{eqnarray}
k_2 = -\frac{1}{\alpha} \left( \left( \alpha' \right)^{\frac{2m}{2m-1}}+1 \right)^{-\frac{1}{2m}}. 
\end{eqnarray}
On the other hand, when $\alpha(u) = u$, we have
\begin{eqnarray}
k_1 = -\frac{1}{2m-1} \left( 1+\left( \beta' \right)^{\frac{2m}{2m-1}} \right)^{-\frac{2m+1}{2m}} (\beta')^{-\frac{2m-2}{2m-1}} \beta'', 
\end{eqnarray}
and
\begin{eqnarray}
k_2 = -\frac{1}{u} \left( 1+\left( \beta' \right)^{\frac{2m}{2m-1}} \right)^{-\frac{1}{2m}} (\beta')^{\frac{1}{2m-1}}. 
\end{eqnarray}

\section{The case where one of principal curvatures is constant}

Let $({\mathbb R}^3, \|\cdot\|)$ be the normed $3$-space as in Section 3. Let $M$ be a rotational surface in $({\mathbb R}^3, \|\cdot\|)$ parametrized as
\[X(u, v) = (\alpha(u)\cos{v}, \alpha(u)\sin{v}, u) \]
where $\alpha > 0$ and $\alpha' \neq 0$. 

In this section, we consider the case where $k_1$ or $k_2$ is constant.

\vspace{2mm}

(i) The case where $k_1 = 0$ identically. In this case, by (3.1), we see that $\alpha(u)$ is a nonconstant linear function, and the surface is a circular cylinder.

\vspace{2mm}

(ii) The case where $k_2 = \mu$ for a constant $\mu$. In this case, we have $\mu < 0$ by (3.2). It suffices to consider the case where $\mu = -1$. Then we have
\[-\frac{1}{\alpha} \left( \left( \alpha' \right)^{\frac{2m}{2m-1}}+1 \right)^{-\frac{1}{2m}} = -1, \]
\[\frac{d\alpha}{du} = \pm \frac{(1-\alpha^{2m})^{\frac{2m-1}{2m}}} {\alpha^{2m-1}}, \]
and
\[u(\alpha) = u_{\pm}(\alpha) := \pm\int \frac{\alpha^{2m-1}} {(1-\alpha^{2m})^{\frac{2m-1}{2m}}} d\alpha = \mp(1-\alpha^{2m})^{\frac{1}{2m}}+c \]
for a constant $c$. It satisfies
\[\alpha^{2m}+(u_{\pm}(\alpha)-c)^{2m} = 1, \]
and the resulting surface is a parallel translation of the unit sphere $S$.

\vspace{2mm}

(iii) The case where $k_1 = \mu$ for a non-zero constant $\mu$. Then by (3.1) we have
\[ \frac{1}{2m-1} \left( \left( \alpha' \right)^{\frac{2m}{2m-1}}+1 \right)^{-\frac{2m+1}{2m}} (\alpha')^{-\frac{2m-2}{2m-1}} \alpha'' = \mu. \]
Multiplying by $-\alpha'$ we have
\[ -\frac{1}{2m-1} \left( \left( \alpha' \right)^{\frac{2m}{2m-1}}+1 \right)^{-\frac{2m+1}{2m}} (\alpha')^{\frac{1}{2m-1}} \alpha'' = -\mu \alpha'. \]
Integrating it we have
\[\left( (\alpha')^{\frac{2m}{2m-1}}+1 \right)^{-\frac{1}{2m}} = c_1-\mu \alpha \ (> 0) \]
for a constant $c_1$. Then
\[\frac{d\alpha}{du} = \pm \frac{\{1-(c_1-\mu \alpha)^{2m}\}^{\frac{2m-1}{2m}}} {(c_1-\mu \alpha)^{2m-1}}, \]
and we get
\[u_{\pm}(\alpha) = \pm \int \frac{(c_1-\mu \alpha)^{2m-1}} {\{1-(c_1-\mu \alpha)^{2m}\}^{\frac{2m-1}{2m}}} d\alpha \]
\[= \pm \frac{1}{\mu} \{1-(c_1-\mu \alpha)^{2m}\}^{\frac{1}{2m}}+c_2 \]
for a constant $c_2$. It suffices to consider the case where $\mu = \pm 1$.

\vspace{2mm}

(iii-1) In the case where $\mu = 1$, we have
\[u_{\pm}(\alpha) = \pm \{1-(c_1-\alpha)^{2m}\}^{\frac{1}{2m}}+c_2, \]
and it satisfies
\begin{eqnarray}
(c_1-\alpha)^{2m}+(u_{\pm}(\alpha)-c_2)^{2m} = 1, \ \ \ \alpha < c_1. 
\end{eqnarray}

\vspace{2mm}

(iii-2) In the case where $\mu = -1$, we have
\[u_{\pm}(\alpha) = \mp \{1-(c_3+\alpha)^{2m}\}^{\frac{1}{2m}}+c_4 \]
for constants $c_3$ and $c_4$. It satisfies
\begin{eqnarray}
(\alpha+c_3)^{2m}+(u_{\pm}(\alpha)-c_4)^{2m} = 1, \ \ \ \alpha > -c_3. 
\end{eqnarray}

If we choose $c_1 > 1$, $c_3 = -c_1$ and $c_4 = c_2$, then (4.1) and (4.2) are smoothly connected. Choosing the orientation suitably, we get a rotational torus with $k_1 = 1$.

\section{Homogeneous linear Weingarten cases}

Let $({\mathbb R}^3, \|\cdot\|)$ be the normed $3$-space as in Section 3. Let $M$ be a rotational surface in $({\mathbb R}^3, \|\cdot\|)$ parametrized as
\[X(u, v) = (\alpha(u)\cos{v}, \alpha(u)\sin{v}, u) \]
where $\alpha > 0$ and $\alpha' \neq 0$. 

In this section, we consider the case where $k_1+\lambda k_2 = 0$ for a non-zero constant $\lambda$. Then by (3.1) and (3.2), we have
\begin{eqnarray}
\frac{1}{2m-1} \left( \left( \alpha' \right)^{\frac{2m}{2m-1}}+1 \right)^{-\frac{2m+1}{2m}} (\alpha')^{-\frac{2m-2}{2m-1}} \alpha'' \nonumber
\end{eqnarray}
\begin{eqnarray}
\hspace{3cm} -\frac{\lambda}{\alpha}\left( \left( \alpha' \right)^{\frac{2m}{2m-1}}+1 \right)^{-\frac{1}{2m}} = 0, 
\end{eqnarray}
and
\[\frac{\alpha''}{(\alpha')^{2}+(\alpha')^{\frac{2m-2}{2m-1}}} = (2m-1)\frac{\lambda}{\alpha}. \]
Multiplying by $2\alpha'$ we have
\[\frac{2\alpha' \alpha''}{(\alpha')^2+(\alpha')^{\frac{2m-2}{2m-1}}} = 2(2m-1)\lambda \cdot \frac{\alpha'}{\alpha}, \]
and
\[2(2m-1)\lambda \log\alpha = \int \frac{((\alpha')^2)'}{(\alpha')^2+(\alpha')^{\frac{2m-2}{2m-1}}} du. \]
Setting
\[(\alpha')^{\frac{2}{2m-1}} =: Z \]
for the right hand side,  we have
\[2\lambda\log\alpha = \int\frac{Z^{m-1}}{Z^{m}+1}dZ =  \frac{1}{m} \log{(Z^{m}+1)}+c_1 \]
\[= \frac{1}{m} \log{((\alpha')^{\frac{2m}{2m-1}}+1)}+c_1 \]
for a constant $c_1$. Then
\[\frac{d\alpha}{du} = \pm \frac{(\alpha^{2m\lambda}-c_{2}^{2m\lambda})^{\frac{2m-1}{2m}}} {c_{2}^{(2m-1)\lambda}} \]
for a positive constant $c_2$, and
\[u(\alpha) = u_{\pm}(\alpha) := \pm \int \frac{c_{2}^{(2m-1)\lambda}} {(\alpha^{2m\lambda}-c_{2}^{2m\lambda})^{\frac{2m-1}{2m}}} d\alpha. \]

\begin{thm}
A rotational surface in $({\mathbb R}^3, \|\cdot\|)$ given by
\[\bar{X}(\alpha, v) = (\alpha\cos{v}, \alpha\sin{v}, u(\alpha)) \]
where $\alpha > 0$ and $u' \neq 0$, satisfies $k_1+\lambda k_2 = 0$ for a non-zero constant $\lambda$, if and only if
\[u(\alpha) = \pm \int \frac{c_{2}^{(2m-1)\lambda}} {(\alpha^{2m\lambda}-c_{2}^{2m\lambda})^{\frac{2m-1}{2m}}} d\alpha \]
for a positive constant $c_2$. 
\end{thm}

(i) The case where $\lambda > 0$. In this case, we have $\alpha > c_2 (> 0)$ and
\[u_{\pm}(\alpha) = \pm \int_{c_2}^{\alpha} \frac{c_{2}^{(2m-1)\lambda}} {(\rho^{2m\lambda}-c_{2}^{2m\lambda})^{\frac{2m-1}{2m}}} d\rho+c_3 \]
for a constant $c_3$. On the denominaor, we set $f(t) := t^{2m\lambda}-c_{2}^{2m\lambda}$. Then
\[f'(t) = 2m\lambda t^{2m\lambda-1}, \ \ \ f(c_2) = 0, \ \ \ f'(c_2) \neq 0. \]
So we have $f(t) = (t-c_2)\tilde{f}(t)$ for some smooth function $\tilde{f}(t)$ such that $\tilde{f}(c_2) \neq 0$. 

Hence, noting that
\[0 < \frac{2m-1}{2m} < 1, \]
we can see that the above integral converses and
\[\lim_{\alpha\rightarrow c_2} \int_{c_2}^{\alpha} \frac{c_{2}^{(2m-1)\lambda}} {(\rho^{2m\lambda}-c_{2}^{2m\lambda})^{\frac{2m-1}{2m}}} d\rho = 0. \]

\vspace{2mm}

(i-1) The case where $(2m-1)\lambda > 1$. In this case, we have
\[\lim_{\alpha\rightarrow\infty} \int_{c_2}^{\alpha} \frac{c_{2}^{(2m-1)\lambda}} {(\rho^{2m\lambda}-c_{2}^{2m\lambda})^{\frac{2m-1}{2m}}} d\rho = d_1 \]
for some positive value $d_1$. Then
\[\lim_{\alpha\rightarrow c_2} u_{\pm}(\alpha) = c_3, \ \ \ \ \lim_{\alpha\rightarrow\infty} u_{\pm}(\alpha) = c_3 \pm d_1, \]
\[\lim_{\alpha\rightarrow c_2} u_{\pm}' (\alpha) = \pm \infty, \ \ \ \ \lim_{\alpha\rightarrow\infty} u_{\pm}' (\alpha) = 0. \]

Let $\alpha_{+}(u)$ be the inverse function of $u_{+}(\alpha)$, which is increasing on $(c_3, c_3+d_1)$ and satisfies
\[\lim_{u\rightarrow c_3} \alpha_{+}(u) = c_2, \ \ \ \ \lim_{u\rightarrow c_3+d_1} \alpha_{+}(u) = \infty, \ \ \ \ \lim_{u\rightarrow c_3} \alpha_{+}' (u) = 0. \]
Let $\alpha_{-}(u)$ be the inverse function of $u_{-}(\alpha)$, which is decreasing on $(c_3-d_1, c_3)$ and satisfies
\[\lim_{u\rightarrow c_3} \alpha_{-}(u) = c_2, \ \ \ \ \lim_{u\rightarrow c_3-d_1} \alpha_{-}(u) = \infty, \ \ \ \ \lim_{u\rightarrow c_3} \alpha_{-}' (u) = 0. \]

We define a function $\hat{\alpha}(u)$ on $(c_3-d_1, c_3+d_1)$ by
\[\hat{\alpha}(u) = \left\{ 
\begin{array}{l}
\alpha_{+}(u), \ \ \ c_3 < u < c_3+d_1 \\[2mm]
\alpha_{-}(u), \ \ \ c_3-d_1 < u < c_3 \\[2mm]
c_2, \ \ \ u = c_3. 
\end{array} \right. \]
It is a $C^1$-function on $(c_3-d_1, c_3+d_1)$ such that
\[\hat{\alpha}' (u) = \left\{ 
\begin{array}{l}
\alpha_{+}' (u), \ \ \ c_3 < u < c_3+d_1 \\[2mm]
\alpha_{-}' (u), \ \ \ c_3-d_1 < u < c_3 \\[2mm]
0, \ \ \ u = c_3. 
\end{array} \right. \]
For $u \in (c_3-d_1, c_3) \cup (c_3, c_3+d_1)$, $\hat{\alpha}(u)$ satisfies the equation (5.1). Then we can see that 
\[\lim_{u\rightarrow c_3} (\hat{\alpha}' (u))^{-\frac{2m-2}{2m-1}} \hat{\alpha}'' (u) = \frac{(2m-1)\lambda}{c_2}, \]
and since $m \geq 2$, 
\[\lim_{u\rightarrow c_3} \hat{\alpha}'' (u) = 0. \]
Hence, $\hat{\alpha}(u)$ is a $C^2$-function on $(c_3-d_1, c_3+d_1)$. 

\begin{thm}
Under the notation above in (i-1), the rotational surface in $({\mathbb R}^3, \|\cdot\|)$ given by
\[\hat{X}(u, v) = (\hat{\alpha}(u)\cos{v}, \hat{\alpha}(u)\sin{v}, u), \ \ \ (u, v) \in (c_3-d_1, c_3+d_1) \times [0, 2\pi] \]
satisfies $k_1+\lambda k_2 = 0$ for a constant $\lambda$ with $(2m-1)\lambda > 1$. 
\end{thm}

(i-2) The case where $0 < (2m-1)\lambda \leq 1$. In this case, we have
\[\lim_{\alpha\rightarrow\infty} \int_{c_2}^{\alpha} \frac{c_{2}^{(2m-1)\lambda}} {(\rho^{2m\lambda}-c_{2}^{2m\lambda})^{\frac{2m-1}{2m}}} d\rho = \infty. \]
Then
\[\lim_{\alpha\rightarrow c_2} u_{\pm}(\alpha) = c_3, \ \ \ \ \lim_{\alpha\rightarrow\infty} u_{\pm}(\alpha) = \pm \infty, \]
\[\lim_{\alpha\rightarrow c_2} u_{\pm}' (\alpha) = \pm \infty, \ \ \ \ \lim_{\alpha\rightarrow\infty} u_{\pm}' (\alpha) = 0. \]

Let $\alpha_{+}(u)$ be the inverse function of $u_{+}(\alpha)$, which is increasing on $(c_3, \infty)$ and satisfies
\[\lim_{u\rightarrow c_3} \alpha_{+}(u) = c_2, \ \ \ \ \lim_{u\rightarrow \infty} \alpha_{+}(u) = \infty, \ \ \ \ \lim_{u\rightarrow c_3} \alpha_{+}' (u) = 0. \]
Let $\alpha_{-}(u)$ be the inverse function of $u_{-}(\alpha)$, which is decreasing on $(-\infty, c_3)$ and satisfies
\[\lim_{u\rightarrow c_3} \alpha_{-}(u) = c_2, \ \ \ \ \lim_{u\rightarrow -\infty} \alpha_{-}(u) = \infty, \ \ \ \ \lim_{u\rightarrow c_3} \alpha_{-}' (u) = 0. \]

We define a function $\hat{\alpha}(u)$ on $(-\infty, \infty)$ by
\[\hat{\alpha}(u) = \left\{ 
\begin{array}{l}
\alpha_{+}(u), \ \ \ c_3 < u < \infty \\[2mm]
\alpha_{-}(u), \ \ \ -\infty < u < c_3 \\[2mm]
c_2, \ \ \ u = c_3. 
\end{array} \right. \]
As in the case (i-1), we can see that $\hat{\alpha}(u)$ is a $C^2$-function on $(-\infty, \infty)$. 

\begin{thm}
Under the notation above in (i-2), the rotational surface in $({\mathbb R}^3, \|\cdot\|)$ given by
\[\hat{X}(u, v) = (\hat{\alpha}(u)\cos{v}, \hat{\alpha}(u)\sin{v}, u), \ \ \ (u, v) \in (-\infty, \infty) \times [0, 2\pi] \]
satisfies $k_1+\lambda k_2 = 0$ for a constant $\lambda$ with $0 < (2m-1)\lambda \leq 1$. 
\end{thm}

(ii) The case where $\lambda < 0$. In this case, we have $0 < \alpha < c_2$ and
\[u_{\pm}(\alpha) = \pm \int_{c_2}^{\alpha} \frac{c_{2}^{(2m-1)\lambda}} {(\rho^{2m\lambda}-c_{2}^{2m\lambda})^{\frac{2m-1}{2m}}} d\rho+c_3 \]
\[= \pm \int_{c_2}^{\alpha} \frac{\rho^{(2m-1)(-\lambda)}} {\{ c_{2}^{2m(-\lambda)}-\rho^{2m(-\lambda)} \}^{\frac{2m-1}{2m}}} d\rho+c_3. \]
As in the case (i), this integral converges and
\[\lim_{\alpha\rightarrow c_2} \int_{c_2}^{\alpha} \frac{\rho^{(2m-1)(-\lambda)}} {\{ c_{2}^{2m(-\lambda)}-\rho^{2m(-\lambda)} \}^{\frac{2m-1}{2m}}} d\rho = 0. \]

Set
\[d_2 := -\lim_{\alpha\rightarrow 0} \int_{c_2}^{\alpha} \frac{\rho^{(2m-1)(-\lambda)}} {\{ c_{2}^{2m(-\lambda)}-\rho^{2m(-\lambda)} \}^{\frac{2m-1}{2m}}} d\rho. \]
Then
\[\lim_{\alpha\rightarrow c_2} u_{\pm}(\alpha) = c_3, \ \ \ \ \lim_{\alpha\rightarrow 0} u_{\pm}(\alpha) = c_3 \mp d_2, \]
\[\lim_{\alpha\rightarrow c_2} u_{\pm}' (\alpha) = \pm\infty, \ \ \ \ \lim_{\alpha\rightarrow 0} u_{\pm}' (\alpha) = 0. \]
Correspondingly, we get a closed rotational surface of $C^1$-class which is homeomorphic to the $2$-sphere. 

As in the case (i), we can see that the surface is of $C^2$-class at points where $\alpha = c_2$. On the $C^2$-smoothness at $\alpha = 0$, we recall that
\[u_{\pm}' (\alpha) = \pm \frac{\alpha^{(2m-1)(-\lambda)}} {\{ c_{2}^{2m(-\lambda)}-\alpha^{2m(-\lambda)} \}^{\frac{2m-1}{2m}}}, \ \ \ \alpha \in (0, c_2). \]
Hence, the limit $\displaystyle{ \lim_{\alpha\rightarrow 0} u_{\pm}'' (\alpha) }$ exists if and only if $(2m-1)(-\lambda) \geq 1$. 

Moreover, using (3.3) and (3.4), we can see that $k_1$ and $k_2$ are extended continuously to $\alpha = 0$ if and only if $(-\lambda) \geq 1$. 

\begin{thm}
There are closed rotational surface such that $k_1+\lambda k_2 = 0$ for a constant $\lambda \leq -1$. They are homeomorphic to the $2$-sphere. 
\end{thm}

\section{Inhomogeneous linear Weingarten cases}

Let $({\mathbb R}^3, \|\cdot\|)$ be the normed $3$-space as in Section 3. Let $M$ be a rotational surface in $({\mathbb R}^3, \|\cdot\|)$ parametrized as
\[X(u, v) = (\alpha(u)\cos{v}, \alpha(u)\sin{v}, u) \]
where $\alpha > 0$ and $\alpha' \neq 0$. 

In this section, we consider the case where $k_1+\lambda k_2 = \mu$ for non-zero constants $\lambda$ and $\mu$. The case where $\lambda = 1$ (the constant mean curvature case) was treated in the previous paper [8]. Then by (3.1) and (3.2), we have
\begin{eqnarray}
\frac{1}{2m-1} \left( \left( \alpha' \right)^{\frac{2m}{2m-1}}+1 \right)^{-\frac{2m+1}{2m}} (\alpha')^{-\frac{2m-2}{2m-1}} \alpha'' \nonumber
\end{eqnarray}
\begin{eqnarray}
\hspace{3cm} -\frac{\lambda}{\alpha}\left( \left( \alpha' \right)^{\frac{2m}{2m-1}}+1 \right)^{-\frac{1}{2m}} = \mu. 
\end{eqnarray}
Multiplying by $-\alpha^{\lambda}\alpha'$ we have
\[\alpha^{\lambda} \left\{ -\frac{1}{2m-1} \left( \left( \alpha' \right)^{\frac{2m}{2m-1}}+1 \right)^{-\frac{2m+1}{2m}} (\alpha')^{\frac{1}{2m-1}} \alpha'' \right\} \]
\[+\lambda\alpha^{\lambda-1}\alpha' \left( \left( \alpha' \right)^{\frac{2m}{2m-1}}+1 \right)^{-\frac{1}{2m}} = -\mu\alpha^{\lambda}\alpha', \]
and
\[\alpha^{\lambda} \left( \left( \alpha' \right)^{\frac{2m}{2m-1}}+1 \right)^{-\frac{1}{2m}} = -\mu \int \alpha^{\lambda}\alpha' du. \]

\subsection{The case where $\lambda = -1$}

In this case, we have
\[\frac{1}{\alpha} \left( \left( \alpha' \right)^{\frac{2m}{2m-1}}+1 \right)^{-\frac{1}{2m}} = -\mu\int\frac{\alpha'}{\alpha}du = -\mu\log{\alpha}+c_1 \ (> 0) \]
for a constant $c_1$. Then
\[\frac{d\alpha}{du} = \pm \frac{\{1-\alpha^{2m}(c_{1}-\mu\log{\alpha})^{2m}\}^{\frac{2m-1}{2m}}} {\alpha^{2m-1}(c_{1}-\mu\log{\alpha})^{2m-1}}, \]
and
\[u(\alpha) = u_{\pm}(\alpha) := \pm \int \frac{\alpha^{2m-1}(c_{1}-\mu\log{\alpha})^{2m-1}} {\{1-\alpha^{2m}(c_{1}-\mu\log{\alpha})^{2m}\}^{\frac{2m-1}{2m}}} d\alpha. \]

\begin{thm}
A rotational surface in $({\mathbb R}^3, \|\cdot\|)$ given by
\[\bar{X}(\alpha, v) = (\alpha\cos{v}, \alpha\sin{v}, u(\alpha)) \]
where $\alpha > 0$ and $u' \neq 0$, satisfies $k_1-k_2 = \mu$ for a non-zero constant $\mu$, if and only if
\[u(\alpha) = \pm \int \frac{\alpha^{2m-1}(c_{1}-\mu\log{\alpha})^{2m-1}} {\{1-\alpha^{2m}(c_{1}-\mu\log{\alpha})^{2m}\}^{\frac{2m-1}{2m}}} d\alpha \]
for a constant $c_1$. Here we need to have
\[c_{1}-\mu\log{\alpha} > 0, \ \ \ \alpha(c_{1}-\mu\log{\alpha}) < 1. \]
\end{thm}

It suffices to consider the case where $\mu = \pm 1$.

\vspace{2mm}

(i) The case where $\mu = 1$. We have
\[u_{\pm}(\alpha) = \pm \int \frac{\alpha^{2m-1}(c_{1}-\log{\alpha})^{2m-1}} {\{1-\alpha^{2m}(c_{1}-\log{\alpha})^{2m}\}^{\frac{2m-1}{2m}}} d\alpha \]
where $0 < \alpha < e^{c_1}$ and $\alpha(c_{1}-\log{\alpha}) < 1$. Set
\[g(t) := t(c_{1}-\log{t}), \ \ \ 0 < t < e^{c_1}. \]
Then
\[g'(t) = c_{1}-1-\log{t}, \]
and we have the following table: 

\vspace{2mm}

\begin{tabular}{|c|c|c|c|c|c|} \hline
$t$ & $(0)$ & $\cdots$ & $e^{c_{1}-1}$ & $\cdots$ & $(e^{c_1})$ \\ \hline 
$g'(t)$ &   & $+$ & $0$ & $-$ &    \\ \hline 
$g(t)$ &   & $\nearrow$ & $e^{c_{1}-1}$ & $\searrow$ &   \\ \hline 
\end{tabular}
\vspace{2mm}
\\
with 
\[\lim_{t\rightarrow 0} g(t) = 0, \ \ \ \ \lim_{t\rightarrow e^{c_1}} g(t) = 0. \]

\vspace{2mm}

(i-1) The case where $c_1 < 1$. In this case, we have $\alpha(c_{1}-\log{\alpha}) < 1$ for $0 < \alpha < e^{c_1}$. So the domain is given by $0 < \alpha < e^{c_1}$, and 
\[u_{\pm}(\alpha) = \pm \int_{e^{c_1}}^{\alpha} \frac{\rho^{2m-1}(c_{1}-\log{\rho})^{2m-1}} {\{1-\rho^{2m}(c_{1}-\log{\rho})^{2m}\}^{\frac{2m-1}{2m}}} d\rho+c_{2}^{\pm}. \]
Set
\[d_1 := -\lim_{\alpha\rightarrow 0} \int_{e^{c_1}}^{\alpha} \frac{\rho^{2m-1}(c_{1}-\log{\rho})^{2m-1}} {\{1-\rho^{2m}(c_{1}-\log{\rho})^{2m}\}^{\frac{2m-1}{2m}}} d\rho. \]
Then we have
\[\lim_{\alpha\rightarrow e^{c_1}} u_{\pm}(\alpha) = c_{2}^{\pm}, \ \ \ \ \lim_{\alpha\rightarrow 0} u_{\pm}(\alpha) = c_{2}^{\pm} \mp d_1, \]
\[\lim_{\alpha\rightarrow e^{c_1}} u_{\pm}' (\alpha) = 0, \ \ \ \ \lim_{\alpha\rightarrow 0} u_{\pm}' (\alpha) = 0. \]

On the $C^2$-smoothness at $\alpha = 0$, we recall that
\[u_{\pm}' (\alpha) = \pm \frac{ \{\alpha(c_{1}-\log{\alpha}) \}^{2m-1} } { [ 1-\{ \alpha(c_{1}-\log{\alpha}) \}^{2m} ]^{\frac{2m-1}{2m}} }, \ \ \ \alpha \in (0, e^{c_1}). \]
Since $m \geq 2$, we can see that $\displaystyle{ \lim_{\alpha\rightarrow 0} u_{\pm}'' (\alpha) }$ exists. However, by (3.3) and (3.4), we find that $k_1$ and $k_2$ are not extended continuously to $\alpha = 0$.

\vspace{2mm}

(i-2) The case where $c_1 = 1$. In this case, the domain is given by $0 < \alpha < 1$ or $1 < \alpha < e$. 

\vspace{1mm}

(i-2-1) On the domain $0 < \alpha < 1$, we have
\[u_{\pm}(\alpha) = \pm \int_{0}^{\alpha} \frac{\rho^{2m-1}(1-\log{\rho})^{2m-1}} {\{1-\rho^{2m}(1-\log{\rho})^{2m}\}^{\frac{2m-1}{2m}}} d\rho+c_{3}^{\pm}. \]
On the demoninator, we set
\[f(t) := 1-t(1-\log{t}). \]
Then
\[f'(t) = \log{t}, \ \ \ f''(t) = \frac{1}{t}, \ \ \ f(1) = 0, \ \ \ f'(1) = 0, \ \ \ f''(1) = 1 \neq 0. \]
So we have $f(t) = (1-t)^{2} \tilde{f}(t)$ for a smooth function $\tilde{f}(t)$ such that $\tilde{f}(1) \neq 0$. 

Hence, noting that $m \geq 2$ and
\[1 < 2\times \frac{2m-1}{2m} < 2, \]
we find that
\[\lim_{\alpha\rightarrow 1} \int_{0}^{\alpha} \frac{\rho^{2m-1}(1-\log{\rho})^{2m-1}} {\{1-\rho^{2m}(1-\log{\rho})^{2m}\}^{\frac{2m-1}{2m}}} d\rho = \infty. \]
Thus we have
\[\lim_{\alpha\rightarrow 0} u_{\pm}(\alpha) = c_{3}^{\pm}, \ \ \ \ \lim_{\alpha\rightarrow 1} u_{\pm}(\alpha) = \pm \infty, \ \ \ \ \lim_{\alpha\rightarrow 0} u_{\pm}'(\alpha) = 0. \]
Correspondingly, we get a rotational surface of $C^1$-class. 

The $C^2$-smoothness at $\alpha = 0$ is analogous to the case (i-1). We can see that $\displaystyle{ \lim_{\alpha\rightarrow 0} u_{\pm}'' (\alpha) }$ exists, but $k_1$ and $k_2$ are not extended continuously to $\alpha = 0$.

\vspace{2mm}

(i-2-2) On the domain $1 < \alpha < e$, we have
\[u_{\pm}(\alpha) = \pm \int_{e}^{\alpha} \frac{\rho^{2m-1}(1-\log{\rho})^{2m-1}} {\{1-\rho^{2m}(1-\log{\rho})^{2m}\}^{\frac{2m-1}{2m}}} d\rho+c_{3}^{\pm}. \]
As in the case (i-2-1), we can see that
\[\lim_{\alpha\rightarrow 1} \int_{e}^{\alpha} \frac{\rho^{2m-1}(1-\log{\rho})^{2m-1}} {\{1-\rho^{2m}(1-\log{\rho})^{2m}\}^{\frac{2m-1}{2m}}} d\rho = -\infty, \]
and
\[\lim_{\alpha\rightarrow e} u_{\pm}(\alpha) = c_{3}^{\pm}, \ \ \ \ \lim_{\alpha\rightarrow 1} u_{\pm}(\alpha) = \mp \infty, \ \ \ \ \lim_{\alpha\rightarrow e} u_{\pm}'(\alpha) = 0. \]

\vspace{1mm}

(i-3) The case where $c_1 > 1$. In this case, there exist two numbers $\alpha_1$ and $\alpha_2$ such that
\[\alpha_{1}(c_{1}-\log{\alpha_1}) = 1, \ \ \ 0 < \alpha_1 < e^{c_{1}-1}, \]
\[\alpha_{2}(c_{1}-\log{\alpha_2}) = 1, \ \ \ e^{c_{1}-1}  < \alpha_2 < e^{c_1}, \]
and the domain is given by $0 < \alpha < \alpha_1$ or $\alpha_2 < \alpha < e^{c_1}$.

\vspace{2mm}

(i-3-1) On the domain $0 < \alpha < \alpha_1$, we have
\[u_{\pm}(\alpha) = \pm \int_{\alpha_1}^{\alpha} \frac{\rho^{2m-1}(c_{1}-\log{\rho})^{2m-1}} {\{1-\rho^{2m}(c_{1}-\log{\rho})^{2m}\}^{\frac{2m-1}{2m}}} d\rho+c_4. \]
On the denominator, we set
\[f(t) := 1-t(c_{1}-\log{t}). \]
Then
\[f'(t) = 1-c_{1}+\log{t}, \ \ \ f(\alpha_1) = 0, \]
\[f'(\alpha_1) = 1-c_{1}+\log{\alpha_1} = 1-\frac{1}{\alpha_1}, \]
where we use that $\alpha_{1}(c_{1}-\log{\alpha_1}) = 1$ for the last equality. Since $c_1 > 1$, we can see that $\alpha_1 \neq 1$, and $f'(\alpha_1) \neq 0$. Thus we have $f(t) = (\alpha_{1}-t)\tilde{f}(t)$ for a smooth function $\tilde{f}(t)$ with $\tilde{f}(\alpha_1) \neq 0$. 

Since
\[0 < \frac{2m-1}{2m} < 1, \]
the above integral converges and
\[\lim_{\alpha\rightarrow \alpha_1} \int_{\alpha_1}^{\alpha} \frac{\rho^{2m-1}(c_{1}-\log{\rho})^{2m-1}} {\{1-\rho^{2m}(c_{1}-\log{\rho})^{2m}\}^{\frac{2m-1}{2m}}} d\rho = 0. \]
Set
\[d_2 := -\lim_{\alpha\rightarrow 0} \int_{\alpha_1}^{\alpha} \frac{\rho^{2m-1}(c_{1}-\log{\rho})^{2m-1}} {\{1-\rho^{2m}(c_{1}-\log{\rho})^{2m}\}^{\frac{2m-1}{2m}}} d\rho. \]
Then
\[\lim_{\alpha\rightarrow \alpha_1} u_{\pm}(\alpha) = c_4, \ \ \ \ \lim_{\alpha\rightarrow 0} u_{\pm}(\alpha) = c_{4} \mp d_2, \]
\[\lim_{\alpha\rightarrow \alpha_1} u_{\pm}' (\alpha) = \pm \infty, \ \ \ \ \lim_{\alpha\rightarrow 0} u_{\pm}' (\alpha) = 0. \]
Correspondingly, we get a rotational surface of $C^1$-class which is homeomorphic to the $2$-sphere. 

Using (6.1), the $C^2$-smoothness at $\alpha = \alpha_1$ can be seen as in Section 5. The $C^2$-smoothness at $\alpha = 0$ is analogous to the case (i-1). We see that $\displaystyle{ \lim_{\alpha\rightarrow 0} u_{\pm}'' (\alpha) }$ exists, but $k_1$ and $k_2$ are not extended continuously to $\alpha = 0$.

\vspace{2mm}

(i-3-2) On the domain $\alpha_2 < \alpha < e^{c_1}$, we have
\[u_{\pm}(\alpha) = \pm \int_{\alpha_2}^{\alpha} \frac{\rho^{2m-1}(c_{1}-\log{\rho})^{2m-1}} {\{1-\rho^{2m}(c_{1}-\log{\rho})^{2m}\}^{\frac{2m-1}{2m}}} d\rho+c_4. \]
As in the case (i-3-1), it converges. Set
\[d_3 := \lim_{\alpha\rightarrow e^{c_1}} \int_{\alpha_2}^{\alpha} \frac{\rho^{2m-1}(c_{1}-\log{\rho})^{2m-1}} {\{1-\rho^{2m}(c_{1}-\log{\rho})^{2m}\}^{\frac{2m-1}{2m}}} d\rho. \]
Then
\[\lim_{\alpha\rightarrow \alpha_2} u_{\pm}(\alpha) = c_4, \ \ \ \ \lim_{\alpha\rightarrow e^{c_1}} u_{\pm}(\alpha) = c_{4}\pm d_3, \]
\[\lim_{\alpha\rightarrow \alpha_2} u_{\pm}' (\alpha) = \pm \infty, \ \ \ \ \lim_{\alpha\rightarrow e^{c_1}} u_{\pm}' (\alpha) = 0. \]
The $C^2$-smoothness at $\alpha = \alpha_2$ can be seen as before.

\vspace{2mm}

(ii) The case where $\mu = -1$. We have
\[u_{\pm}(\alpha) = \pm \int \frac{\alpha^{2m-1}(c_{1}^{\ast}+\log{\alpha})^{2m-1}} {\{1-\alpha^{2m}(c_{1}^{\ast}+\log{\alpha})^{2m}\}^{\frac{2m-1}{2m}}} d\alpha \]
where $\alpha > e^{-c_{1}^{\ast}}$ and $\alpha(c_{1}^{\ast}+\log{\alpha}) < 1$. Set
\[g(t) := t(c_{1}^{\ast}+\log{t}), \ \ \ t > e^{-c_{1}^{\ast}}. \]
Then, for $t > e^{-c_{1}^{\ast}}$, we have
\[g'(t) = c_{1}^{\ast}+1+\log{t} >1, \ \ \ \lim_{t\rightarrow e^{-c_{1}^{\ast}}} g(t) = 0, \ \ \ \lim_{t\rightarrow\infty} g(t) = \infty. \]
So there exists a unique number $\alpha_3$ such that
\[\alpha_{3}(c_{1}^{\ast}+\log{\alpha_3}) = 1, \ \ \ \alpha_3 > e^{-c_{1}^{\ast}}, \]
and the domain is given by $e^{-c_{1}^{\ast}} < \alpha < \alpha_3$. Thus we have
\[u_{\pm}(\alpha) = \pm \int_{\alpha_3}^{\alpha} \frac{\rho^{2m-1}(c_{1}^{\ast}+\log{\rho})^{2m-1}} {\{1-\rho^{2m}(c_{1}^{\ast}+\log{\rho})^{2m}\}^{\frac{2m-1}{2m}}} d\rho+c_5. \]

On the denominator, we set
\[f(t) := 1-t(c_{1}^{\ast}+\log{t}). \]
Then
\[f'(t) = -(c_{1}^{\ast}+1+\log{t}), \ \ \ f(\alpha_3) = 0, \]
\[f'(\alpha_3) = -(c_{1}^{\ast}+1+\log{\alpha_{3}}) = -\left( 1+\frac{1}{\alpha_3} \right) < -1, \]
where we use that $\alpha_{3}(c_{1}^{\ast}+\log{\alpha_3}) = 1$ for the last equality. Thus we have $f(t) = (\alpha_{3}-t)\tilde{f}(t)$ for a smooth function $\tilde{f}(t)$ with $\tilde{f}(\alpha_3) \neq 0$. 

Hence, the above integral converges. Set
\[d_4 := -\lim_{\alpha\rightarrow e^{-c_{1}^{\ast}}} \int_{\alpha_3}^{\alpha} \frac{\rho^{2m-1}(c_{1}^{\ast}+\log{\rho})^{2m-1}} {\{1-\rho^{2m}(c_{1}^{\ast}+\log{\rho})^{2m}\}^{\frac{2m-1}{2m}}} d\rho. \]
Then
\[\lim_{\alpha\rightarrow \alpha_3} u_{\pm}(\alpha) = c_5, \ \ \ \ \lim_{\alpha\rightarrow e^{-c_{1}^{\ast}}} u_{\pm}(\alpha) = c_{5} \mp d_4, \]
\[\lim_{\alpha\rightarrow \alpha_3} u_{\pm}' (\alpha) = \pm \infty, \ \ \ \ \lim_{\alpha\rightarrow e^{-c_{1}^{\ast}}} u_{\pm}' (\alpha) = 0. \]
The $C^2$-smoothness at $\alpha = \alpha_3$ can be seen as before.

\subsection{Connecting argument for $\lambda = -1$}

Let $u_{1\pm}$ denote the $u_{\pm}$ in the case (i) of Subsection 6.1, and $u_{2\pm}$ denote the $u_{\pm}$ in the case (ii) of Subsection 6.1.

\vspace{2mm}

(C1) Connecting (i-1) and (ii) of Subsection 6.1.

\vspace{2mm}

(C1-1) If we choose 
\[c_{1}^{\ast} = -c_1, \ \ \ c_{2}^{+} = c_{5}-d_4, \ \ \ c_{2}^{-} = c_{5}+d_4, \]
then by connecting the surfaces in the cases (i-1) and (ii), we get a closed rotational surface of $C^1$-class which is homeomorphic to the $2$-sphere. Notice that $\alpha_3$ depends on $c_1$. 

Noting the orientation, by (3.3) and (3.4) for $\beta = u_{1+}(\alpha)$ and $\beta = u_{2+}(\alpha)$, we have
\[ -\frac{1}{2m-1} \left( 1+\left( u_{1+}' \right)^{\frac{2m}{2m-1}} \right)^{-\frac{2m+1}{2m}} (u_{1+}')^{-\frac{2m-2}{2m-1}} u_{1+}'' \hspace{2cm} \]
\[ \hspace{2cm} +\frac{1}{\alpha} \left( 1+\left( u_{1+}' \right)^{\frac{2m}{2m-1}} \right)^{-\frac{1}{2m}} (u_{1+}')^{\frac{1}{2m-1}} = 1, \ \ \ \ \alpha \in (0, e^{c_1}), \]
and
\[ -\frac{1}{2m-1} \left( 1+\left( u_{2+}' \right)^{\frac{2m}{2m-1}} \right)^{-\frac{2m+1}{2m}} (u_{2+}')^{-\frac{2m-2}{2m-1}} u_{2+}'' \hspace{2cm} \]
\[ \hspace{2cm} +\frac{1}{\alpha} \left( 1+\left( u_{2+}' \right)^{\frac{2m}{2m-1}} \right)^{-\frac{1}{2m}} (u_{2+}')^{\frac{1}{2m-1}} = -1, \ \ \ \ \alpha \in (e^{c_1}, \alpha_3). \]
Thus we find that
\[\lim_{\alpha\rightarrow e^{c_1}} (u_{1+}' (\alpha))^{-\frac{2m-2}{2m-1}} u_{1+}'' (\alpha) = -(2m-1) \]
and
\[\lim_{\alpha\rightarrow e^{c_1}} (u_{2+}' (\alpha))^{-\frac{2m-2}{2m-1}} u_{2+}'' (\alpha) = 2m-1. \]
Since $m \geq 2$, we have
\[\lim_{\alpha\rightarrow e^{c_1}} u_{1+}'' (\alpha) = \lim_{\alpha\rightarrow e^{c_1}} u_{2+}'' (\alpha) = 0. \]
Hence, $u_{1+}(\alpha)$ and $u_{2+}(\alpha)$ are connected $C^2$-smoothly at $\alpha = e^{c_1}$, but $k_1$ is not continuous at $\alpha = e^{c_1}$. The connectedness for $u_{1-}(\alpha)$ and $u_{2-}(\alpha)$ is analogous.

\vspace{2mm}

(C1-2) On the other hand, choose 
\[c_{1}^{\ast} = -c_1, \ \ \ c_{2}^{+} = c_{5}+d_4, \ \ \ c_{2}^{-} = c_{5}-d_4. \]
If $d_1 \neq d_4$, then we have a rotational surface of $C^1$-class. 

In this case, noting the difference of the orientation, we have
\[ -\frac{1}{2m-1} \left( 1+\left( u_{1+}' \right)^{\frac{2m}{2m-1}} \right)^{-\frac{2m+1}{2m}} (u_{1+}')^{-\frac{2m-2}{2m-1}} u_{1+}'' \hspace{2cm} \]
\[ \hspace{2cm} +\frac{1}{\alpha} \left( 1+\left( u_{1+}' \right)^{\frac{2m}{2m-1}} \right)^{-\frac{1}{2m}} (u_{1+}')^{\frac{1}{2m-1}} = 1, \ \ \ \ \alpha \in (0, e^{c_1}), \]
and
\[ -\frac{1}{2m-1} \left( 1+\left( u_{2-}' \right)^{\frac{2m}{2m-1}} \right)^{-\frac{2m+1}{2m}} (u_{2-}')^{-\frac{2m-2}{2m-1}} u_{2-}'' \hspace{2cm} \]
\[ \hspace{2cm} +\frac{1}{\alpha} \left( 1+\left( u_{2-}' \right)^{\frac{2m}{2m-1}} \right)^{-\frac{1}{2m}} (u_{2-}')^{\frac{1}{2m-1}} = 1, \ \ \ \ \alpha \in (e^{c_1}, \alpha_3). \]
We can see that $u_{1+}(\alpha)$ and $u_{2-}(\alpha)$ are connected $C^2$-smoothly at $\alpha = e^{c_1}$, and $k_1$ is continuous at $\alpha = e^{c_1}$ in this case. The connectedness for $u_{1-}(\alpha)$ and $u_{2+}(\alpha)$ is analogous. 

This surface is of $C^2$-class, but as mentioned in (i-1) of Subsection 6.1, $k_1$ and $k_2$ are not extended continuously to $\alpha = 0$. 

\vspace{2mm}

{\bf Remark.} If $d_1 = d_4$, then it has a singularity at a point where $\alpha = 0$.

\vspace{2mm}

(C2) Connecting (i-2-2) and (ii) of Subsection 6.1. If we choose 
\[c_{1}^{\ast} = -1, \ \ \ c_{3}^{+} = c_{5}+d_4, \ \ \ c_{3}^{-} = c_{5}-d_4, \]
then we have a rotational surface of $C^1$-class. 

As in the case (C1-2), $u_{1+}(\alpha)$ and $u_{2-}(\alpha)$ are connected $C^2$-smoothly at $\alpha = e$, and $k_1$ is continuous at $\alpha = e$. The connectedness for $u_{1-}(\alpha)$ and $u_{2+}(\alpha)$ is analogous.

\vspace{2mm}

(C3) Connecting (i-3-2) and (ii) of Subsection 6.1. If we choose
\[c_{1}^{\ast} = -c_1, \ \ \ c_4 = c_{5}+d_{4}-d_3, \]
then we have a rotational surface of $C^1$-class. Notice that $\alpha_2$ and $\alpha_3$ depend on $c_1$. 

As in the case (C1-2) and (C2), we can see that $u_{1+}(\alpha)$ and $u_{2-}(\alpha)$ are connected $C^2$-smoothly at $\alpha = e^{c_1}$, and $k_1$ is continuous at $\alpha = e^{c_1}$. Analogously, we can see that $u_{1-}(\alpha)$ and $u_{2+}(\alpha)$ have the same derivatives at $\alpha = e^{c_1}$. 

Let $\Gamma$ be the $C^2$-curve given by connecting the graphs of $u_{1-}$, $u_{1+}$, $u_{2-}$ and $u_{2+}$. Since $\Gamma$ has the same derivatives at the end points, it can be extended periodically as a $C^2$-curve $\Gamma^{\ast}$. As a result, we obtain a periodic rotational surface which satisfies $k_1-k_2 = 1$ for a suitable choice of orientation.

\vspace{2mm}

{\bf Remark.} If $d_3 = d_4$, then we may have a rotational torus such that $k_1-k_2$ is a non-zero constant.

\subsection{The case where $\lambda \neq -1$}

In this case, we have
\[\alpha^{\lambda} \left( \left( \alpha' \right)^{\frac{2m}{2m-1}}+1 \right)^{-\frac{1}{2m}} = -\mu\int \alpha^{\lambda}\alpha' du =  -\frac{\mu}{\lambda+1} \alpha^{\lambda+1}+c_1 \ (> 0) \]
for a constant $c_1$. Then
\[\frac{d\alpha}{du} = \pm \frac{ [ (\lambda+1)^{2m}\alpha^{2m\lambda}-\{ c_1(\lambda+1)-\mu\alpha^{\lambda+1} \}^{2m} ]^{\frac{2m-1}{2m}}} {\{ c_1(\lambda+1)-\mu\alpha^{\lambda+1} \}^{2m-1}} \]
and
\[u(\alpha) = \pm \int \frac{\{ c_1(\lambda+1)-\mu\alpha^{\lambda+1} \}^{2m-1}} {[ (\lambda+1)^{2m}\alpha^{2m\lambda}-\{ c_1(\lambda+1)-\mu\alpha^{\lambda+1} \}^{2m} ]^{\frac{2m-1}{2m}}} d\alpha. \]

\begin{thm}
A rotational surface in $({\mathbb R}^3, \|\cdot\|)$ given by
\[\bar{X}(\alpha, v) = (\alpha\cos{v}, \alpha\sin{v}, u(\alpha)) \]
where $\alpha > 0$ and $u' \neq 0$, satisfies $k_1+\lambda k_2 = \mu$ for non-zero constants $\lambda \neq -1$ and $\mu$, if and only if
\[u(\alpha) = \pm \int \frac{\{ c_1(\lambda+1)-\mu\alpha^{\lambda+1} \}^{2m-1}} {[ (\lambda+1)^{2m}\alpha^{2m\lambda}-\{ c_1(\lambda+1)-\mu\alpha^{\lambda+1} \}^{2m} ]^{\frac{2m-1}{2m}}} d\alpha \]
for a constant $c_1$. Here we need to have
\[0 < c_{1}-\frac{\mu}{\lambda+1} \alpha^{\lambda+1} < \alpha^{\lambda}. \]
\end{thm}

We divide the cases as $\lambda > 0$, $-1 < \lambda < 0$ or $\lambda < -1$. It suffices to consider the case where $\mu = \pm 1$.

\vspace{2mm}

(i) The case where $\lambda > 0$ and $\mu = 1$. In this case, we have
\[u_{\pm}(\alpha) = \pm \int \frac{\{ c_1(\lambda+1)-\alpha^{\lambda+1} \}^{2m-1}} {[ (\lambda+1)^{2m}\alpha^{2m\lambda}-\{ c_1(\lambda+1)-\alpha^{\lambda+1} \}^{2m} ]^{\frac{2m-1}{2m}}} d\alpha \]
and
\begin{eqnarray}
0 < c_{1}(\lambda+1)-\alpha^{\lambda+1} < (\lambda+1)\alpha^{\lambda}. 
\end{eqnarray}
From the first inequality of (6.2), we need to have $c_1 > 0$, and
\[0 < \alpha < \{c_1 (\lambda+1)\}^{\frac{1}{\lambda+1}} =: \alpha_{4}. \]
For the second inequality of (6.2), set
\[f(t) := (\lambda+1) t^{\lambda}-c_1(\lambda+1)+t^{\lambda+1}, \ \ \ 0 < t < \alpha_{4}. \]
Since it is strictly increasing and
\[\lim_{t\rightarrow 0} f(t) = -c_1(\lambda+1) < 0, \ \ \ \ \lim_{t\rightarrow \alpha_4} f(t) = (\lambda+1)\alpha_{4}^{\lambda} > 0, \]
there exists a unique number $\alpha_5 \in (0, \alpha_4)$ such that $f(\alpha_5) = 0$. 

Thus the domain is given by $\alpha_{5} < \alpha < \alpha_{4}$, and
\[u_{\pm}(\alpha) = \pm \int_{\alpha_5}^{\alpha} \frac{\{ c_1(\lambda+1)-\rho^{\lambda+1} \}^{2m-1}} {[ (\lambda+1)^{2m}\rho^{2m\lambda}-\{ c_1(\lambda+1)-\rho^{\lambda+1} \}^{2m} ]^{\frac{2m-1}{2m}}} d\rho+c_2. \]
On the denominator, since $f(\alpha_5) = 0$ and $f'(\alpha_5) > 0$, we have $f(t) = (t-\alpha_5)\tilde{f}(t)$ for a smooth function $\tilde{f}(t)$ with $\tilde{f}(\alpha_5) \neq 0$. So the above integral converges. 

Set
\[d_1 := \lim_{\alpha\rightarrow \alpha_4} \int_{\alpha_5}^{\alpha} \frac{\{ c_1(\lambda+1)-\rho^{\lambda+1} \}^{2m-1}} {[ (\lambda+1)^{2m}\rho^{2m\lambda}-\{ c_1(\lambda+1)-\rho^{\lambda+1} \}^{2m} ]^{\frac{2m-1}{2m}}} d\rho. \]
Then we have
\[\lim_{\alpha\rightarrow \alpha_5} u_{\pm}(\alpha) = c_2, \ \ \ \ \lim_{\alpha\rightarrow \alpha_4} u_{\pm}(\alpha) = c_2\pm d_1, \]
\[\lim_{\alpha\rightarrow \alpha_5} u_{\pm}' (\alpha) = \pm\infty, \ \ \ \ \lim_{\alpha\rightarrow \alpha_4} u_{\pm}' (\alpha) = 0. \]

\vspace{2mm}

(ii) The case where $\lambda > 0$ and $\mu = -1$. We have
\[u_{\pm}(\alpha) = \pm \int \frac{\{ c_{1}^{\ast} (\lambda+1)+\alpha^{\lambda+1} \}^{2m-1}} {[ (\lambda+1)^{2m}\alpha^{2m\lambda}-\{ c_{1}^{\ast} (\lambda+1)+\alpha^{\lambda+1} \}^{2m} ]^{\frac{2m-1}{2m}}} d\alpha \]
and
\begin{eqnarray}
0 < c_{1}^{\ast} (\lambda+1)+\alpha^{\lambda+1} < (\lambda+1)\alpha^{\lambda}. 
\end{eqnarray}

\vspace{1mm}

(ii-1) The case where $c_{1}^{\ast} = 0$. We have
\[u_{\pm}(\alpha) = \pm \int \frac{ \alpha^{(2m-1)(\lambda+1)} } { \{ (\lambda+1)^{2m}\alpha^{2m\lambda}-\alpha^{2m(\lambda+1)} \}^{\frac{2m-1}{2m}} } d\alpha \]
\[= \pm \int \frac{\alpha^{2m-1}} {\{(\lambda+1)^{2m}-\alpha^{2m}\}^{\frac{2m-1}{2m}}} d\alpha = \mp \{(\lambda+1)^{2m}-\alpha^{2m}\}^{\frac{1}{2m}}+c_3. \]
It satisfies
\[\alpha^{2m}+(u_{\pm}(\alpha)-c_3)^{2m} = (\lambda+1)^{2m}, \]
and the surface is homothetic to the unit sphere $S$.

\vspace{2mm}

(ii-2) The case where $c_{1}^{\ast} > 0$. In this case, the first inequality of (6.3) holds obviously. For the second inequality of (6.3), set
\[f(t) := (\lambda+1)t^{\lambda}-c_{1}^{\ast} (\lambda+1)-t^{\lambda+1}, \ \ \ t > 0. \]
Then
\[f'(t) = (\lambda+1)(\lambda-t)t^{\lambda-1} \]
and
\[\lim_{t\rightarrow 0} f(t) = -c_{1}^{\ast} (\lambda+1) < 0, \ \ \ \ \lim_{t\rightarrow\infty} f(t) = -\infty. \]
So we need to have 
\[f(\lambda) = \lambda^{\lambda}-c_{1}^{\ast} (\lambda+1) > 0, \]
that is, 
\[0 < c_{1}^{\ast} < \frac{\lambda^{\lambda}}{\lambda+1}. \]
Then there are two numbers $\alpha_6$ and $\alpha_7$ such that 
\[f(\alpha_6) = f(\alpha_7) = 0, \ \ \ \ 0 <\alpha_6 < \lambda < \alpha_7. \]

Thus the domain is given by $\alpha_6 < \alpha < \alpha_7$, and
\[u_{\pm}(\alpha) = \pm \int_{\alpha_6}^{\alpha} \frac{\{ c_{1}^{\ast} (\lambda+1)+\rho^{\lambda+1} \}^{2m-1}} {[ (\lambda+1)^{2m}\rho^{2m\lambda}-\{ c_{1}^{\ast} (\lambda+1)+\rho^{\lambda+1} \}^{2m} ]^{\frac{2m-1}{2m}}} d\rho+c_3. \]
On the denominator, we have
\[f(\alpha_6) = f(\alpha_7) = 0, \ \ \ f'(\alpha_6) \neq 0, \ \ \ f'(\alpha_7) \neq 0. \]
So the above integral converges as $\alpha$ tends to $\alpha_6$ and $\alpha_7$. 

Set
\[d_2 := \lim_{\alpha\rightarrow \alpha_7} \int_{\alpha_6}^{\alpha} \frac{\{ c_{1}^{\ast} (\lambda+1)+\rho^{\lambda+1} \}^{2m-1}} {[ (\lambda+1)^{2m}\rho^{2m\lambda}-\{ c_{1}^{\ast} (\lambda+1)+\rho^{\lambda+1} \}^{2m} ]^{\frac{2m-1}{2m}}} d\rho. \]
Then
\[\lim_{\alpha\rightarrow \alpha_6} u_{\pm}(\alpha) = c_3, \ \ \ \ \lim_{\alpha\rightarrow \alpha_7} u_{\pm}(\alpha) = c_3 \pm d_2, \]
\[\lim_{\alpha\rightarrow \alpha_6} u_{\pm}' (\alpha) = \pm\infty, \ \ \ \ \lim_{\alpha\rightarrow \alpha_7} u_{\pm}' (\alpha) = \pm\infty. \]

Let $\Gamma$ be the $C^2$-curve given by connecting $u_{+}$ and $u_{-}$. Since $\Gamma$ has the same derivatives at the end points, it can be extended periodically as a $C^2$-curve $\Gamma^{\ast}$. Thus we get a periodic rotational surface.

\vspace{2mm}

(ii-3) The case where $c_{1}^{\ast} < 0$. From the first inequality of (6.3), we have
\[\alpha > \{(-c_{1}^{\ast})(\lambda+1)\}^{\frac{1}{\lambda+1}} =: \alpha_8. \]
For the second inequality of (6.3), set
\[f(t) := (\lambda+1)t^{\lambda}-c_{1}^{\ast} (\lambda+1)-t^{\lambda+1}, \ \ \ t > \alpha_8. \]
Then
\[f'(t) = (\lambda+1)(\lambda-t)t^{\lambda-1} \]
and
\[\lim_{t\rightarrow \alpha_8} f(t) = (\lambda+1)\alpha_{8}^{\lambda} > 0, \ \ \ \ \lim_{t\rightarrow\infty} f(t) = -\infty. \]
So there is a unique number $\alpha_9 \in (\alpha_8, \infty)$ such that $f(\alpha_9) = 0$. 

Thus the domain is given by $\alpha_8 < \alpha < \alpha_9$, and
\[u_{\pm}(\alpha) = \pm \int_{\alpha_9}^{\alpha} \frac{\{ c_{1}^{\ast}(\lambda+1)+\rho^{\lambda+1} \}^{2m-1}} {[ (\lambda+1)^{2m}\rho^{2m\lambda}-\{ c_{1}^{\ast}(\lambda+1)+\rho^{\lambda+1} \}^{2m} ]^{\frac{2m-1}{2m}}} d\rho+c_3. \]
On the denominator, we can see that $f(\alpha_9) = 0$ and $f'(\alpha_9) \neq 0$, so that the above integral converges. 

Set
\[d_3 := -\lim_{\alpha\rightarrow \alpha_8} \int_{\alpha_9}^{\alpha} \frac{\{ c_{1}^{\ast} (\lambda+1)+\rho^{\lambda+1} \}^{2m-1}} {[ (\lambda+1)^{2m}\rho^{2m\lambda}-\{ c_{1}^{\ast} (\lambda+1)+\rho^{\lambda+1} \}^{2m} ]^{\frac{2m-1}{2m}}} d\rho. \]
Then
\[\lim_{\alpha\rightarrow \alpha_9} u_{\pm}(\alpha) = c_3, \ \ \ \ \lim_{\alpha\rightarrow \alpha_8} u_{\pm}(\alpha) = c_3 \mp d_3, \]
\[\lim_{\alpha\rightarrow \alpha_9} u_{\pm}' (\alpha) = \pm\infty, \ \ \ \ \lim_{\alpha\rightarrow \alpha_8} u_{\pm}' (\alpha) = 0. \]

\vspace{2mm}

(iii) The case where $-1 < \lambda < 0$ and $\mu = 1$. We have
\[u_{\pm}(\alpha) = \pm \int \frac{\{ c_1(\lambda+1)-\alpha^{\lambda+1} \}^{2m-1}} {[ (\lambda+1)^{2m}\alpha^{2m\lambda}-\{ c_1(\lambda+1)-\alpha^{\lambda+1} \}^{2m} ]^{\frac{2m-1}{2m}}} d\alpha \]
\[= \pm \int \frac{ \alpha^{(2m-1)(-\lambda)} \{ c_1(\lambda+1)-\alpha^{\lambda+1} \}^{2m-1}} {[ (\lambda+1)^{2m}-\alpha^{2m(-\lambda)} \{ c_1(\lambda+1)-\alpha^{\lambda+1} \}^{2m} ]^{\frac{2m-1}{2m}}} d\alpha, \]
and
\begin{eqnarray}
c_1(\lambda+1)-\alpha^{\lambda+1} > 0, \ \ \ \alpha^{(-\lambda)} \{c_{1}(\lambda+1)-\alpha^{\lambda+1}\} < \lambda+1. 
\end{eqnarray}
From the first inequality of (6.4), we have $c_1 > 0$ and
\[0 < \alpha < \{c_1(\lambda+1)\}^{\frac{1}{\lambda+1}} =: \alpha_{10}. \]
For the second inequality of (6.4), set
\[f(t) := \lambda+1-c_{1}(\lambda+1) t^{(-\lambda)}+t, \ \ \ 0 < t < \alpha_{10}. \]
Then
\[f'(t) = 1-c_1(-\lambda)(\lambda+1)t^{-(\lambda+1)}. \]
Set
\[\alpha_{11} := \{c_1(-\lambda)(\lambda+1)\}^{\frac{1}{\lambda+1}} < \alpha_{10}, \]
where we use that $0 <-\lambda < 1$ for the inequality. 

So we have the following table: 

\vspace{2mm}

\begin{tabular}{|c|c|c|c|c|c|} \hline
$t$ & $(0)$ & $\cdots$ & $\alpha_{11}$ & $\cdots$ & $(\alpha_{10})$ \\ \hline 
$f'(t)$ &   & $-$ & $0$ & $+$ &    \\ \hline 
$f(t)$ &   & $\searrow$ & $   $ & $\nearrow$ &   \\ \hline 
\end{tabular}
\vspace{2mm}
\\
with
\[\lim_{t\rightarrow 0} f(t) = \lim_{t\rightarrow \alpha_{10}} f(t) = \lambda+1 > 0, \]
and
\[f(\alpha_{11}) = (\lambda+1)\left[ 1-(-\lambda)^{\frac{(-\lambda)}{\lambda+1}} \{c_1(\lambda+1)\}^{\frac{1}{\lambda+1}} \right]. \]

\vspace{2mm}

(iii-1) The case where $0 < c_1 < \{(\lambda+1)(-\lambda)^{(-\lambda)}\}^{-1}$. In this case, we have $f(\alpha_{11}) > 0$ and the domain is given by $0 < \alpha < \alpha_{10}$. Then
\[u_{\pm}(\alpha) = \pm \int_{0}^{\alpha} \frac{ \rho^{(2m-1)(-\lambda)} \{ c_1(\lambda+1)-\rho^{\lambda+1} \}^{2m-1}} {[ (\lambda+1)^{2m}-\rho^{2m(-\lambda)} \{ c_1(\lambda+1)-\rho^{\lambda+1} \}^{2m} ]^{\frac{2m-1}{2m}}} d\rho+c_{4}^{\pm}. \]
Set
\[d_4 := \lim_{\alpha\rightarrow \alpha_{10}} \int_{0}^{\alpha} \frac{ \rho^{(2m-1)(-\lambda)} \{ c_1(\lambda+1)-\rho^{\lambda+1} \}^{2m-1}} {[ (\lambda+1)^{2m}-\rho^{2m(-\lambda)} \{ c_1(\lambda+1)-\rho^{\lambda+1} \}^{2m} ]^{\frac{2m-1}{2m}}} d\rho. \]
Then
\[\lim_{\alpha\rightarrow 0} u_{\pm}(\alpha) = c_{4}^{\pm}, \ \ \ \ \lim_{\alpha\rightarrow \alpha_{10}} u_{\pm}(\alpha) = c_{4}^{\pm} \pm d_4, \]
\[\lim_{\alpha\rightarrow 0} u_{\pm}' (\alpha) = 0, \ \ \ \ \lim_{\alpha\rightarrow \alpha_{10}} u_{\pm}' (\alpha) = 0. \]

On the $C^2$-smoothness at $\alpha = 0$, we note that
\[u_{\pm}' (\alpha) = \pm\frac{ \{ c_1(\lambda+1)\alpha^{(-\lambda)}-\alpha \}^{2m-1}} {[ (\lambda+1)^{2m}-\{ c_1(\lambda+1)\alpha^{(-\lambda)}-\alpha \}^{2m} ]^{\frac{2m-1}{2m}}}. \]
Hence, $\displaystyle{ \lim_{\alpha\rightarrow 0} u_{\pm}'' (\alpha) }$ exists if and only if $(2m-1)(-\lambda) \geq 1$. 

Howover, by (3.3), (3.4) and that $-1 < \lambda < 0$, we find that $k_1$ and $k_2$ are not extended continuously to $\alpha = 0$.

\vspace{2mm}

(iii-2) The case where $c_1 = \{(\lambda+1)(-\lambda)^{(-\lambda)}\}^{-1}$. In this case, we have $f(\alpha_{11}) = 0$ and the domain is given by $0 < \alpha < \alpha_{11}$ or $\alpha_{11} < \alpha < \alpha_{10}$.

\vspace{2mm}

(iii-2-1) On the domain $0 < \alpha < \alpha_{11}$, we have
\[u_{\pm}(\alpha) = \pm \int_{0}^{\alpha} \frac{ \rho^{(2m-1)(-\lambda)} \{ c_1(\lambda+1)-\rho^{\lambda+1} \}^{2m-1}} {[ (\lambda+1)^{2m}-\rho^{2m(-\lambda)} \{ c_1(\lambda+1)-\rho^{\lambda+1} \}^{2m} ]^{\frac{2m-1}{2m}}} d\rho+c_{4}^{\pm}. \]
On the denominator, since $f(\alpha_{11}) = f'(\alpha_{11}) = 0$ and $f'' (\alpha_{11}) \neq 0$, we have $f(t) = (\alpha_{11}-t)^{2}\tilde{f}(t)$ for a smooth function $\tilde{f}(t)$ with $\tilde{f}(\alpha_{11}) \neq 0$. 

Thus, noting that $m \geq 2$, we have
\[\lim_{\alpha\rightarrow \alpha_{11}} \int_{0}^{\alpha} \frac{ \rho^{(2m-1)(-\lambda)} \{ c_1(\lambda+1)-\rho^{\lambda+1} \}^{2m-1}} {[ (\lambda+1)^{2m}-\rho^{2m(-\lambda)} \{ c_1(\lambda+1)-\rho^{\lambda+1} \}^{2m} ]^{\frac{2m-1}{2m}}} d\rho = \infty, \]
and
\[\lim_{\alpha\rightarrow 0} u_{\pm}(\alpha) = c_{4}^{\pm}, \ \ \ \ \lim_{\alpha\rightarrow \alpha_{11}} u_{\pm}(\alpha) = \pm \infty, \ \ \ \ \lim_{\alpha\rightarrow 0} u_{\pm}' (\alpha) = 0. \]
Correspondingly, we get a rotational surface of $C^1$-class. 

As in the case (iii-1), $\displaystyle{ \lim_{\alpha\rightarrow 0} u_{\pm}'' (\alpha) }$ exists if and only if $(2m-1)(-\lambda) \geq 1$, but $k_1$ and $k_2$ are not extended continuously to $\alpha = 0$.

\vspace{2mm}

(iii-2-2) On the domain $\alpha_{11} < \alpha < \alpha_{10}$, we have
\[u_{\pm}(\alpha) = \pm \int_{\alpha_{10}}^{\alpha} \frac{ \rho^{(2m-1)(-\lambda)} \{ c_1(\lambda+1)-\rho^{\lambda+1} \}^{2m-1}} {[ (\lambda+1)^{2m}-\rho^{2m(-\lambda)} \{ c_1(\lambda+1)-\rho^{\lambda+1} \}^{2m} ]^{\frac{2m-1}{2m}}} d\rho+c_{4}^{\pm}. \]
As in the case (iii-2-1), we can see that
\[\lim_{\alpha\rightarrow \alpha_{10}} u_{\pm}(\alpha) = c_{4}^{\pm}, \ \ \ \ \lim_{\alpha\rightarrow \alpha_{11}} u_{\pm}(\alpha) = \mp \infty, \ \ \ \ \lim_{\alpha\rightarrow \alpha_{10}} u_{\pm}' (\alpha) = 0. \]

\vspace{2mm}

(iii-3) The case where $c_1 > \{(\lambda+1)(-\lambda)^{(-\lambda)}\}^{-1}$. In this case, we have $f(\alpha_{11}) < 0$. There are two numbers $\alpha_{12}$ and $\alpha_{13}$ such that
\[f(\alpha_{12}) = f(\alpha_{13}) = 0, \ \ \ 0 < \alpha_{12} < \alpha_{11} < \alpha_{13} < \alpha_{10}, \]
and the domain is given by $0 < \alpha < \alpha_{12}$ or $\alpha_{13} < \alpha < \alpha_{10}$.

\vspace{2mm}

(iii-3-1) On the domain $0 < \alpha < \alpha_{12}$, we have
\[u_{\pm}(\alpha) = \pm \int_{\alpha_{12}}^{\alpha} \frac{ \rho^{(2m-1)(-\lambda)} \{ c_1(\lambda+1)-\rho^{\lambda+1} \}^{2m-1}} {[ (\lambda+1)^{2m}-\rho^{2m(-\lambda)} \{ c_1(\lambda+1)-\rho^{\lambda+1} \}^{2m} ]^{\frac{2m-1}{2m}}} d\rho+c_{5}. \]
On the denominator, since $f(\alpha_{12}) = 0$ and $f' (\alpha_{12}) \neq 0$, the above integral converges. 

Set
\[d_5 := -\lim_{\alpha\rightarrow 0} \int_{\alpha_{12}}^{\alpha} \frac{ \rho^{(2m-1)(-\lambda)} \{ c_1(\lambda+1)-\rho^{\lambda+1} \}^{2m-1}} {[ (\lambda+1)^{2m}-\rho^{2m(-\lambda)} \{ c_1(\lambda+1)-\rho^{\lambda+1} \}^{2m} ]^{\frac{2m-1}{2m}}} d\rho. \]
Then
\[\lim_{\alpha\rightarrow \alpha_{12}} u_{\pm}(\alpha) = c_{5}, \ \ \ \ \lim_{\alpha\rightarrow 0} u_{\pm}(\alpha) = c_{5} \mp d_5, \]
\[\lim_{\alpha\rightarrow \alpha_{12}} u_{\pm}' (\alpha) = \pm \infty, \ \ \ \ \lim_{\alpha\rightarrow 0} u_{\pm}' (\alpha) = 0. \]
Correspondingly, we get a closed rotational surface of $C^1$-class which is homeomorphic to the $2$-sphere. 

As in the case (iii-1), $\displaystyle{ \lim_{\alpha\rightarrow 0} u_{\pm}'' (\alpha) }$ exists if and only if $(2m-1)(-\lambda) \geq 1$, but $k_1$ and $k_2$ are not extended continuously to $\alpha = 0$.

\vspace{2mm}

(iii-3-2) On the domain $\alpha_{13} < \alpha < \alpha_{10}$, we have
\[u_{\pm}(\alpha) = \pm \int_{\alpha_{13}}^{\alpha} \frac{ \rho^{(2m-1)(-\lambda)} \{ c_1(\lambda+1)-\rho^{\lambda+1} \}^{2m-1}} {[ (\lambda+1)^{2m}-\rho^{2m(-\lambda)} \{ c_1(\lambda+1)-\rho^{\lambda+1} \}^{2m} ]^{\frac{2m-1}{2m}}} d\rho+c_5. \]
As in the case (iii-3-1), it converges. Set
\[d_6 := \lim_{\alpha\rightarrow \alpha_{10}} \int_{\alpha_{13}}^{\alpha} \frac{ \rho^{(2m-1)(-\lambda)} \{ c_1(\lambda+1)-\rho^{\lambda+1} \}^{2m-1}} {[ (\lambda+1)^{2m}-\rho^{2m(-\lambda)} \{ c_1(\lambda+1)-\rho^{\lambda+1} \}^{2m} ]^{\frac{2m-1}{2m}}} d\rho. \]
Then
\[\lim_{\alpha\rightarrow \alpha_{13}} u_{\pm}(\alpha) = c_{5}, \ \ \ \ \lim_{\alpha\rightarrow \alpha_{10}} u_{\pm}(\alpha) = c_{5} \pm d_6, \]
\[\lim_{\alpha\rightarrow \alpha_{13}} u_{\pm}' (\alpha) = \pm \infty, \ \ \ \ \lim_{\alpha\rightarrow \alpha_{10}} u_{\pm}' (\alpha) = 0. \]

\vspace{2mm}

(iv) The case where $-1 < \lambda < 0$ and $\mu = -1$. We have
\[u_{\pm}(\alpha) = \pm \int \frac{\{ c_{1}^{\ast} (\lambda+1)+\alpha^{\lambda+1} \}^{2m-1}} {[ (\lambda+1)^{2m}\alpha^{2m\lambda}-\{ c_{1}^{\ast} (\lambda+1)+\alpha^{\lambda+1} \}^{2m} ]^{\frac{2m-1}{2m}}} d\alpha \]
\[= \pm \int \frac{ \alpha^{(2m-1)(-\lambda)} \{ c_{1}^{\ast} (\lambda+1)+\alpha^{\lambda+1} \}^{2m-1}} {[ (\lambda+1)^{2m}-\alpha^{2m(-\lambda)} \{ c_{1}^{\ast} (\lambda+1)+\alpha^{\lambda+1} \}^{2m} ]^{\frac{2m-1}{2m}}} d\alpha, \]
and
\begin{eqnarray}
c_{1}^{\ast} (\lambda+1)+\alpha^{\lambda+1} > 0, \ \ \ \alpha^{(-\lambda)} \{c_{1}^{\ast} (\lambda+1)+\alpha^{\lambda+1}\} < \lambda+1. 
\end{eqnarray}

\vspace{1mm}

(iv-1) The case where $c_{1}^{\ast} = 0$. As in the case (ii-1), the surface is homothetic to the unit sphere $S$.

\vspace{2mm}

(iv-2) The case where $c_{1}^{\ast} > 0$. The first inequality of (6.5) holds obviously. For the second inequality of (6.5), set
\[f(t) := \lambda+1-c_{1}^{\ast} (\lambda+1) t^{(-\lambda)}-t, \ \ \ t > 0. \]
Since it is strictly decreasing and
\[\lim_{t\rightarrow 0} f(t) = \lambda+1 > 0, \ \ \ \ \lim_{t\rightarrow\infty} f(t) = -\infty, \]
there is a unique number $\alpha_{14} > 0$ such that $f(\alpha_{14}) = 0$. 

The domain is given by $0 < \alpha < \alpha_{14}$, and
\[u_{\pm}(\alpha) = \pm \int_{\alpha_{14}}^{\alpha} \frac{ \rho^{(2m-1)(-\lambda)} \{ c_{1}^{\ast} (\lambda+1)+\rho^{\lambda+1} \}^{2m-1}} {[ (\lambda+1)^{2m}-\rho^{2m(-\lambda)} \{ c_{1}^{\ast} (\lambda+1)+\rho^{\lambda+1} \}^{2m} ]^{\frac{2m-1}{2m}}} d\rho+c_6. \]
On the denominator, since $f(\alpha_{14}) = 0$ and $f'(\alpha_{14}) < 0$, the above integral converges. 

Set
\[d_7 := -\lim_{\alpha\rightarrow 0} \int_{\alpha_{14}}^{\alpha} \frac{ \rho^{(2m-1)(-\lambda)} \{ c_{1}^{\ast} (\lambda+1)+\rho^{\lambda+1} \}^{2m-1}} {[ (\lambda+1)^{2m}-\rho^{2m(-\lambda)} \{ c_{1}^{\ast} (\lambda+1)+\rho^{\lambda+1} \}^{2m} ]^{\frac{2m-1}{2m}}} d\rho. \]
Then
\[\lim_{\alpha\rightarrow \alpha_{14}} u_{\pm}(\alpha) = c_{6}, \ \ \ \ \lim_{\alpha\rightarrow 0} u_{\pm}(\alpha) = c_{6} \mp d_7, \]
\[\lim_{\alpha\rightarrow \alpha_{14}} u_{\pm}' (\alpha) = \pm \infty, \ \ \ \ \lim_{\alpha\rightarrow 0} u_{\pm}' (\alpha) = 0. \]
Thus we have a closed rotational surface of $C^1$-class which is homeomorphic to the $2$-sphere. 

As in the case (iii-1), $\displaystyle{ \lim_{\alpha\rightarrow 0} u_{\pm}'' (\alpha) }$ exists if and only if $(2m-1)(-\lambda) \geq 1$, but $k_1$ and $k_2$ are not extended continuously to $\alpha = 0$.

\vspace{2mm}

(iv-3) The case where $c_{1}^{\ast} < 0$. From the first inequality of (6.5), we have
\[\alpha > \{(-c_{1}^{\ast} )(\lambda+1)\}^{\frac{1}{\lambda+1}} =: \alpha_{15}. \]
For the second inequality of (6.5), set
\[f(t) := \lambda+1-c_{1}^{\ast} (\lambda+1) t^{(-\lambda)}-t, \ \ \ t >\alpha_{15}. \]
Since
\[f'(t) = (-c_{1}^{\ast})(-\lambda)(\lambda+1)t^{-(\lambda+1)}-1 < -1, \]
and
\[\lim_{t\rightarrow \alpha_{15}} f(t) =  \lambda+1 > 0, \ \ \ \ \lim_{t\rightarrow\infty} f(t) = -\infty, \]
there is a unique number $\alpha_{16} \in (\alpha_{15}, \infty)$ such that $f(\alpha_{16}) = 0$. 

The domain is given by $\alpha_{15} < \alpha < \alpha_{16}$, and
\[u_{\pm}(\alpha) = \pm \int_{\alpha_{16}}^{\alpha} \frac{ \rho^{(2m-1)(-\lambda)} \{ c_{1}^{\ast} (\lambda+1)+\rho^{\lambda+1} \}^{2m-1}} {[ (\lambda+1)^{2m}-\rho^{2m(-\lambda)} \{ c_{1}^{\ast} (\lambda+1)+\rho^{\lambda+1} \}^{2m} ]^{\frac{2m-1}{2m}}} d\rho+c_6. \]
On the denominator, since $f(\alpha_{16}) = 0$ and $f'(\alpha_{16}) <-1$, the above integral converges. 

Set
\[d_8 := -\lim_{\alpha\rightarrow \alpha_{15}} \int_{\alpha_{16}}^{\alpha} \frac{ \rho^{(2m-1)(-\lambda)} \{ c_{1}^{\ast} (\lambda+1)+\rho^{\lambda+1} \}^{2m-1}} {[ (\lambda+1)^{2m}-\rho^{2m(-\lambda)} \{ c_{1}^{\ast} (\lambda+1)+\rho^{\lambda+1} \}^{2m} ]^{\frac{2m-1}{2m}}}d\rho. \]
Then
\[\lim_{\alpha\rightarrow \alpha_{16}} u_{\pm}(\alpha) = c_{6}, \ \ \ \ \lim_{\alpha\rightarrow \alpha_{15}} u_{\pm}(\alpha) = c_{6} \mp d_8, \]
\[\lim_{\alpha\rightarrow \alpha_{16}} u_{\pm}' (\alpha) = \pm \infty, \ \ \ \ \lim_{\alpha\rightarrow \alpha_{15}} u_{\pm}' (\alpha) = 0. \]

\vspace{2mm}

(v) The case where $\lambda < -1$ and $\mu = 1$. Set $\omega := -(\lambda+1) > 0$. Then
\[u_{\pm}(\alpha) = \pm \int \frac{(c_{1}\omega+\alpha^{-\omega})^{2m-1}} {\{ \omega^{2m}\alpha^{-2m(\omega+1)}-(c_{1} \omega+\alpha^{-\omega})^{2m} \}^{\frac{2m-1}{2m}}} d\alpha \]
\[= \pm \int \frac{\alpha^{2m-1}(c_{1}\omega\alpha^{\omega}+1)^{2m-1}} {\left\{ \omega^{2m}-\alpha^{2m}(c_{1}\omega\alpha^{\omega}+1)^{2m} \right\}^{\frac{2m-1}{2m}}} d\alpha, \]
where
\begin{eqnarray}
c_{1}\omega\alpha^{\omega}+1 > 0, \ \ \ \ \alpha(c_{1}\omega\alpha^{\omega}+1) < \omega. 
\end{eqnarray}

\vspace{1mm}

(v-1) The case where $c_1 = 0$. The surface is homothetic to the unit sphere $S$.

\vspace{2mm}

(v-2) The case where $c_1 > 0$. The first inequality of (6.6) holds obviously. For the second inequality of (6.6), set
\[f(t) := \omega-t(c_{1}\omega t^{\omega}+1), \ \ \ t > 0. \]
Since
\[f'(t) = -c_{1}\omega(\omega+1)t^{\omega}-1 < -1 \ \ \mbox{for} \ t > 0, \]
and
\[\lim_{t\rightarrow 0} f(t) = \omega > 0, \ \ \ \ \lim_{t\rightarrow\infty} f(t) = -\infty, \]
there is a unique number $\alpha_{17} > 0$ such that $f(\alpha_{17}) = 0$. 

The domain is given by $0 < \alpha < \alpha_{17}$, and
\[u_{\pm}(\alpha) = \pm \int_{\alpha_{17}}^{\alpha} \frac{\rho^{2m-1}(c_{1}\omega\rho^{\omega}+1)^{2m-1}} {\left\{ \omega^{2m}-\rho^{2m}(c_{1}\omega\rho^{\omega}+1)^{2m} \right\}^{\frac{2m-1}{2m}}} d\rho+c_{7}. \]
On the denominator, we have $f(\alpha_{17}) = 0$ and $f'(\alpha_{17}) < -1$, so that the above integral converges. 

Set
\[d_9 := -\lim_{\alpha\rightarrow 0} \int_{\alpha_{17}}^{\alpha} \frac{\rho^{2m-1}(c_{1}\omega\rho^{\omega}+1)^{2m-1}} {\left\{ \omega^{2m}-\rho^{2m}(c_{1}\omega\rho^{\omega}+1)^{2m} \right\}^{\frac{2m-1}{2m}}} d\rho. \]
Then
\[\lim_{\alpha\rightarrow \alpha_{17}} u_{\pm}(\alpha) = c_{7}, \ \ \ \lim_{\alpha\rightarrow 0} u_{\pm}(\alpha) = c_{7} \mp d_9, \]
\[\lim_{\alpha\rightarrow \alpha_{17}} u_{\pm}' (\alpha) = \pm \infty, \ \ \ \lim_{\alpha\rightarrow 0} u_{\pm}' (\alpha) = 0. \]
Thus we have a closed rotational surface of $C^1$-class which is homeomorphic to the $2$-sphere. 

On the $C^2$-smoothness at $\alpha = 0$, we note that
\[u_{\pm}' (\alpha) = \pm \frac{ (c_{1}\omega\alpha^{\omega+1}+\alpha)^{2m-1}} {\{ \omega^{2m}-(c_{1}\omega\alpha^{\omega+1}+\alpha)^{2m} \}^{\frac{2m-1}{2m}}}. \]
In this case, we find that $\displaystyle{ \lim_{\alpha\rightarrow 0} u_{\pm}'' (\alpha) }$ exists, and by (3.3) and (3.4), $k_1$ and $k_2$ are extended continuously to $\alpha = 0$.

\vspace{2mm}

(v-3) The case where $c_1 < 0$. From the first ineaulity of (6.6), we have
\[0 < \alpha < \frac{1}{ \{(-c_1)\omega \}^{1/\omega}} =: \alpha_{18}. \]
For the second inequality of (6.6), set
\[f(t) := \omega-t(c_{1}\omega t^{\omega}+1), \ \ \ 0 < t < \alpha_{18}. \]
Then
\[f'(t) = (-c_{1})\omega(\omega+1)t^{\omega}-1. \]
Set
\[\alpha_{19} := \frac{1}{ \{(-c_1)\omega(\omega+1) \}^{1/\omega}} < \alpha_{18}. \]
We have the following table: 

\vspace{2mm}

\begin{tabular}{|c|c|c|c|c|c|} \hline
$t$ & $(0)$ & $\cdots$ & $\alpha_{19}$ & $\cdots$ & $(\alpha_{18})$ \\ \hline 
$f'(t)$ &   & $-$ & $0$ & $+$ &    \\ \hline 
$f(t)$ &   & $\searrow$ & $   $ & $\nearrow$ &   \\ \hline 
\end{tabular}
\vspace{2mm}
\\
with
\[\lim_{\alpha\rightarrow 0} f(t) = \lim_{\alpha\rightarrow \alpha_{18}} f(t) = \omega \]
and
\[f(\alpha_{19}) = \omega-\frac{\omega}{\omega+1} \cdot \frac{1}{ \{(-c_1)\omega(\omega+1) \}^{1/\omega}}. \]

\vspace{2mm}

(v-3-1) The case where $c_1 < -\frac{1}{\omega(\omega+1)^{\omega+1}}$. In this case, we have $f(\alpha_{19}) > 0$ and the domain is given by $0 < \alpha < \alpha_{18}$. Then
\[u_{\pm}(\alpha) = \pm \int_{0}^{\alpha} \frac{\rho^{2m-1}(c_{1}\omega\rho^{\omega}+1)^{2m-1}} {\left\{ \omega^{2m}-\rho^{2m}(c_{1}\omega\rho^{\omega}+1)^{2m} \right\}^{\frac{2m-1}{2m}}} d\rho+c_{8}^{\pm}. \]
Set
\[d_{10} := \lim_{\alpha\rightarrow \alpha_{18}} \int_{0}^{\alpha} \frac{\rho^{2m-1}(c_{1}\omega\rho^{\omega}+1)^{2m-1}} {\left\{ \omega^{2m}-\rho^{2m}(c_{1}\omega\rho^{\omega}+1)^{2m} \right\}^{\frac{2m-1}{2m}}} d\rho. \]
Then
\[\lim_{\alpha\rightarrow 0} u_{\pm}(\alpha) = c_{8}^{\pm}, \ \ \ \ \lim_{\alpha\rightarrow \alpha_{18}} u_{\pm}(\alpha) = c_{8}^{\pm} \pm d_{10}, \]
\[\lim_{\alpha\rightarrow 0} u_{\pm}' (\alpha) = 0, \ \ \ \ \lim_{\alpha\rightarrow \alpha_{18}} u_{\pm}' (\alpha) = 0. \]

As in the case (v-2), $\displaystyle{ \lim_{\alpha\rightarrow 0} u_{\pm}'' (\alpha) }$ exists, and $k_1$ and $k_2$ are extended continuously to $\alpha = 0$.

\vspace{2mm}

(v-3-2) The case where $c_1 = -\frac{1}{\omega(\omega+1)^{\omega+1}}$. In this case, we have $f(\alpha_{19}) = 0$ and the domain is given by $0 < \alpha < \alpha_{19}$ or $\alpha_{19} < \alpha < \alpha_{18}$. 

\vspace{2mm}

(v-3-2-1) On the domain $0 < \alpha < \alpha_{19}$, we have
\[u_{\pm}(\alpha) = \pm \int_{0}^{\alpha} \frac{\rho^{2m-1}(c_{1}\omega\rho^{\omega}+1)^{2m-1}} {\left\{ \omega^{2m}-\rho^{2m}(c_{1}\omega\rho^{\omega}+1)^{2m} \right\}^{\frac{2m-1}{2m}}} d\rho+c_{8}^{\pm}. \]
On the denominator, since $f(\alpha_{19}) = f' (\alpha_{19}) = 0$, $f'' (\alpha_{19}) \neq 0$ and $m \geq 2$, we have
\[\lim_{\alpha\rightarrow \alpha_{19}} \int_{0}^{\alpha} \frac{\rho^{2m-1}(c_{1}\omega\rho^{\omega}+1)^{2m-1}} {\left\{ \omega^{2m}-\rho^{2m}(c_{1}\omega\rho^{\omega}+1)^{2m} \right\}^{\frac{2m-1}{2m}}} d\rho = \infty. \]
Then
\[\lim_{\alpha\rightarrow 0} u_{\pm}(\alpha) = c_{8}^{\pm}, \ \ \ \ \lim_{\alpha\rightarrow \alpha_{19}} u_{\pm}(\alpha) = \pm \infty, \ \ \ \ \lim_{\alpha\rightarrow 0} u_{\pm}' (\alpha) = 0. \]
Thus we have a rotational surface of $C^1$-class. 

As in the case (v-2), $\displaystyle{ \lim_{\alpha\rightarrow 0} u_{\pm}'' (\alpha) }$ exists, and $k_1$ and $k_2$ are extended continuously to $\alpha = 0$.

\vspace{2mm}

(v-3-2-2) On the domain $\alpha_{19} < \alpha < \alpha_{18}$, we have
\[u_{\pm}(\alpha) = \pm \int_{\alpha_{18}}^{\alpha} \frac{\rho^{2m-1}(c_{1}\omega\rho^{\omega}+1)^{2m-1}} {\left\{ \omega^{2m}-\rho^{2m}(c_{1}\omega\rho^{\omega}+1)^{2m} \right\}^{\frac{2m-1}{2m}}} d\rho+c_{8}^{\pm}. \]
As in the case (v-3-2-1), we can see that
\[\lim_{\alpha\rightarrow \alpha_{18}} u_{\pm}(\alpha) = c_{8}^{\pm}, \ \ \ \ \lim_{\alpha\rightarrow \alpha_{19}} u_{\pm}(\alpha) = \mp \infty, \ \ \ \ \lim_{\alpha\rightarrow \alpha_{18}} u_{\pm}' (\alpha) = 0. \]

\vspace{2mm}

(v-3-3) The case where $-\frac{1}{\omega(\omega+1)^{\omega+1}} < c_1 < 0$. In this case, we have $f(\alpha_{19}) < 0$. So there are two numbers $\alpha_{20}$ and $\alpha_{21}$ such that
\[f(\alpha_{20}) = f(\alpha_{21}) = 0, \ \ \ 0 < \alpha_{20} < \alpha_{19} < \alpha_{21} < \alpha_{18}, \]
and the domain is given by $0 < \alpha < \alpha_{20}$ or $\alpha_{21} < \alpha < \alpha_{18}$. 

\vspace{2mm}

(v-3-3-1) On the domain $0 < \alpha < \alpha_{20}$, we have
\[u_{\pm}(\alpha) = \pm \int_{\alpha_{20}}^{\alpha} \frac{\rho^{2m-1}(c_{1}\omega\rho^{\omega}+1)^{2m-1}} {\left\{ \omega^{2m}-\rho^{2m}(c_{1}\omega\rho^{\omega}+1)^{2m} \right\}^{\frac{2m-1}{2m}}} d\rho+c_9. \]
On the denominator, since $f(\alpha_{20}) = 0$ and $f' (\alpha_{20}) \neq 0$, the above integral converges. 

Set
\[d_{11} := -\lim_{\alpha\rightarrow 0} \int_{\alpha_{20}}^{\alpha} \frac{\rho^{2m-1}(c_{1}\omega\rho^{\omega}+1)^{2m-1}} {\left\{ \omega^{2m}-\rho^{2m}(c_{1}\omega\rho^{\omega}+1)^{2m} \right\}^{\frac{2m-1}{2m}}} d\rho. \]
Then
\[\lim_{\alpha\rightarrow \alpha_{20}} u_{\pm}(\alpha) = c_9, \ \ \ \ \lim_{\alpha\rightarrow 0} u_{\pm}(\alpha) = c_{9} \mp d_{11}, \]
\[\lim_{\alpha\rightarrow \alpha_{20}} u_{\pm}' (\alpha) = \pm \infty, \ \ \ \ \lim_{\alpha\rightarrow 0} u_{\pm}' (\alpha) = 0. \]
Thus we have a closed rotational surface of $C^1$-class which is homeomorphic to the $2$-sphere. 

As in the case (v-2), $\displaystyle{ \lim_{\alpha\rightarrow 0} u_{\pm}'' (\alpha) }$ exists, and $k_1$ and $k_2$ are extended continuously to $\alpha = 0$.

\vspace{2mm}

(v-3-3-2) On the doamin $\alpha_{21} < \alpha < \alpha_{18}$, we have
\[u_{\pm}(\alpha) = \pm \int_{\alpha_{21}}^{\alpha} \frac{\rho^{2m-1}(c_{1}\omega\rho^{\omega}+1)^{2m-1}} {\left\{ \omega^{2m}-\rho^{2m}(c_{1}\omega\rho^{\omega}+1)^{2m} \right\}^{\frac{2m-1}{2m}}} d\rho+c_9. \]
As in the case (v-3-3-1), it converges. Set
\[d_{12} := \lim_{\alpha\rightarrow \alpha_{18}} \int_{\alpha_{21}}^{\alpha} \frac{\rho^{2m-1}(c_{1}\omega\rho^{\omega}+1)^{2m-1}} {\left\{ \omega^{2m}-\rho^{2m}(c_{1}\omega\rho^{\omega}+1)^{2m} \right\}^{\frac{2m-1}{2m}}} d\rho. \]
Then
\[\lim_{\alpha\rightarrow \alpha_{21}} u_{\pm}(\alpha) = c_9, \ \ \ \ \lim_{\alpha\rightarrow \alpha_{18}} u_{\pm}(\alpha) = c_{9} \pm d_{12}, \]
\[\lim_{\alpha\rightarrow \alpha_{21}} u_{\pm}' (\alpha) = \pm \infty, \ \ \ \ \lim_{\alpha\rightarrow \alpha_{18}} u_{\pm}' (\alpha) = 0. \]

\vspace{2mm}

(vi) The case where $\lambda < -1$ and $\mu = -1$. Setting $\omega := -(\lambda+1) > 0$, we have
\[u_{\pm}(\alpha) = \pm \int \frac{(c_{1}^{\ast} \omega-\alpha^{-\omega})^{2m-1}} {\{ \omega^{2m}\alpha^{-2m(\omega+1)}-(c_{1}^{\ast} \omega-\alpha^{-\omega})^{2m} \}^{\frac{2m-1}{2m}}} d\alpha \]
\[= \pm \int \frac{\alpha^{2m-1}(c_{1}^{\ast} \omega\alpha^{\omega}-1)^{2m-1}} {\left\{ \omega^{2m}-\alpha^{2m}(c_{1}^{\ast} \omega\alpha^{\omega}-1)^{2m} \right\}^{\frac{2m-1}{2m}}} d\alpha, \]
where
\begin{eqnarray}
c_{1}^{\ast} \omega\alpha^{\omega}-1 > 0, \ \ \ \ \alpha(c_{1}^{\ast} \omega\alpha^{\omega}-1) < \omega. 
\end{eqnarray}
From the first inequality of (6.7), we have $c_{1}^{\ast} > 0$ and
\[\alpha > \frac{1}{(c_{1}^{\ast} \omega)^{1/\omega}} =: \alpha_{22}. \]
For the second inequality of (6.7), set
\[f(t) := \omega-t(c_{1}^{\ast} \omega t^{\omega}-1), \ \ \ t > \alpha_{22}. \]
Since
\[f'(t) = 1-c_{1}^{\ast} \omega(\omega+1)t^{\omega} < -\omega \ \ \mbox{for} \ t > \alpha_{22}, \]
and
\[\lim_{t\rightarrow \alpha_{22}} f(t) = \omega > 0, \ \ \ \lim_{t\rightarrow\infty} f(t) = -\infty, \]
there is a unique number $\alpha_{23} \in(\alpha_{22}, \infty)$ such that $f(\alpha_{23}) = 0$. 

The domain is given by $\alpha_{22} < \alpha < \alpha_{23}$, and
\[u_{\pm}(\alpha) = \pm \int_{\alpha_{23}}^{\alpha} \frac{\rho^{2m-1}(c_{1}^{\ast} \omega \rho^{\omega}-1)^{2m-1}} {\left\{ \omega^{2m}-\rho^{2m}(c_{1}^{\ast} \omega \rho^{\omega}-1)^{2m} \right\}^{\frac{2m-1}{2m}}} d\rho+c_{10}. \]
Since $f(\alpha_{23}) = 0$ and $f' (\alpha_{23}) \neq 0$, this integral converges. 

Set
\[d_{13} := -\lim_{\alpha\rightarrow \alpha_{22}} \int_{\alpha_{23}}^{\alpha} \frac{\rho^{2m-1}(c_{1}^{\ast} \omega \rho^{\omega}-1)^{2m-1}} {\left\{ \omega^{2m}-\rho^{2m}(c_{1}^{\ast} \omega \rho^{\omega}-1)^{2m} \right\}^{\frac{2m-1}{2m}}} d\rho. \]
Then
\[\lim_{\alpha\rightarrow \alpha_{23}} u_{\pm}(\alpha) = c_{10}, \ \ \ \ \lim_{\alpha\rightarrow \alpha_{22}} u_{\pm}(\alpha) = c_{10} \mp d_{13}, \]
\[\lim_{\alpha\rightarrow \alpha_{23}} u_{\pm}' (\alpha) = \pm \infty, \ \ \ \ \lim_{\alpha\rightarrow \alpha_{22}} u_{\pm}' (\alpha) = 0. \]

\subsection{Connecting argument for $\lambda \neq -1$}

We use the notation $u_{1\pm}$, $u_{2\pm}$,..., $u_{6\pm}$ to denote the $u_{\pm}$ in the case (i), (ii),..., (vi) of Subsection 6.3, respectievly.

\vspace{2mm}

(C4) Connecting (i) and (ii-3) of Subsection 6.3. Choose
\[c_{1}^{\ast} = -c_1, \ \ \ c_2 = c_{3}+d_{3}-d_1. \]
Then $u_{1+}$ and $u_{2-}$ are connected at $\alpha = \alpha_4 = \alpha_8$. As in Subsection 6.2, we can see that they are connected $C^2$-smoothly, and also $k_1$ and $k_2$ are continuous. 

Let $\Gamma$ be the $C^2$-curve given by connecting the graphs of $u_{1-}$, $u_{1+}$, $u_{2-}$ and $u_{2+}$. Since $\Gamma$ has the same derivatives at the end points, it can be extended periodically as a $C^2$-curve $\Gamma^{\ast}$, and we get a periodic rotational surface. If $d_1 = d_3$, then we may have a rotational torus such that $k_1+\lambda k_2$ is a non-zero constant, where $\lambda$ is a positive constant.

\vspace{2mm}

(C5) Connecting (iii-1) and (iv-3) of Subsection 6.3. Choose
\[c_{1}^{\ast} = -c_1, \ \ \ c_{4}^{+} = c_{6}+d_{8}-d_4, \ \ \ c_{4}^{-} = c_{6}-d_{8}+d_4. \]
If $(2m-1)(-\lambda) > 1$ and $d_4 \neq d_8$, we get a rotational surface of $C^2$-class, but $k_1$ and $k_2$ are not extended continuously to $\alpha = 0$.

\vspace{2mm}

(C6) Connecting (iii-2-2) and (iv-3) of Subsection 6.3. Choose
\[c_{1}^{\ast} = -c_1, \ \ \ c_{4}^{+} = c_{6}+d_8, \ \ \ c_{4}^{-} = c_{6}-d_8. \]
Then we have a rotational surface of $C^2$-class.

\vspace{2mm}

(C7) Connecting (iii-3-2) and (iv-3) of Subsection 6.3. Choose
\[c_{1}^{\ast} = -c_1, \ \ \ c_5 = c_{6}+d_{8}-d_6. \]
In this case, we can get a periodic rotational surface. If $d_6 = d_8$, then we may have a rotational torus such that $k_1+\lambda k_2$ is a non-zero constant, where $\lambda$ is a constant with $-1 < \lambda < 0$.

\vspace{2mm}

(C8) Connecting (v-3-1) and (vi) of Subsection 6.3. Choose
\[c_{1}^{\ast} = -c_1, \ \ \ c_{8}^{+} = c_{10}+d_{13}-d_{10}, \ \ \ c_{8}^{-} = c_{10}-d_{13}+d_{10}. \]
If $d_{10} \neq d_{13}$, then we have a rotational surface of $C^2$-class, and $k_1$ and $k_2$ are extended continuously to $\alpha = 0$.

\vspace{2mm}

(C9) Connecting (v-3-2-2) and (vi) of Subsection 6.3. Choose
\[c_{1}^{\ast} = -c_1, \ \ \ c_{8}^{+} = c_{10}+d_{13}, \ \ \ c_{8}^{-} = c_{10}-d_{13}. \]
Then we have a rotational surface of $C^2$-class.

\vspace{2mm}

(C10) Connecting (v-3-3-2) and (vi) of Subsection 6.3. Choose
\[c_{1}^{\ast} = -c_1, \ \ \ c_9 = c_{10}+d_{13}-d_{12}. \]
In this case, we have a periodic rotational surface. If $d_{12} = d_{13}$, then we may have a rotational torus such that $k_1+\lambda k_2$ is a non-zero constant, where $\lambda$ is a constant with $\lambda < -1$.

\vspace{3mm}

Graduate School of Science and Technology 

Hirosaki University 

Hirosaki 036-8561, Japan 

E-mail: sakaki@hirosaki-u.ac.jp 

\end{document}